\documentclass{amsart}

\usepackage{color,graphics}
\usepackage[matrix,arrow]{xy}
\usepackage[driverfallback=hypertex]{hyperref}

\numberwithin{equation}{section}

\theoremstyle{plain}
\newtheorem{thm}{Theorem}[section]
\newtheorem{assumption}[thm]{Assumption}

\newtheorem{theorem}[thm]{Theorem}

\newtheorem{corollary}[thm]{Corollary}
\newtheorem{lemma}[thm]{Lemma}

\newtheorem{prop}[thm]{Proposition}
\newtheorem{proposition}[thm]{Proposition}

\newtheorem{definition}[thm]{Definition}

\newtheorem{remark}[thm]{Remark}

\newtheorem{example}[thm]{Example}

\newcommand{\R}{\mathbb{R}}
\newcommand{\Z}{\mathbb{Z}}
\newcommand{\Q}{\mathbb{Q}}
\newcommand{\C}{\mathbb{C}}

\newcommand{\half}{{\textstyle\frac{1}{2}}}

\newcommand{\iso}{\cong}
\newcommand{\smooth}{C^\infty}

\newcommand{\A}{\mathcal{A}}
\newcommand{\TT}{\mathcal{T}}
\newcommand{\K}{\mathbb{K}}
\newcommand{\CF}{\mathit{CF}}
\newcommand{\HF}{\mathit{HF}}
\newcommand{\HH}{\mathit{HH}}
\newcommand{\CC}{\mathit{CC}}
\newcommand{\F}{\mathcal{F}}
\newcommand{\SH}{\mathit{SH}}
\newcommand{\Ext}{\mathit{Ext}}
\newcommand{\scrL}{\mathcal{L}}
\newcommand{\id}{\mathrm{id}}
\newcommand{\eq}{\mathrm{eq}}
\newcommand{\scrK}{\mathcal{K}}
\newcommand{\scrM}{\mathcal{M}}

\newcommand{\Hom}{\mathit{Hom}}
\newcommand{\checked}{}

\newcommand{\scrR}{\mathcal{R}}

\begin{document}
\title[$q$-intersection numbers]{Symplectic cohomology and $q$-intersection numbers}
\author{Paul Seidel and Jake P. Solomon}
\date{This revision: March 2012}

\begin{abstract}
Given a symplectic cohomology class of degree $1$, we define the notion of an ``equivariant'' Lagrangian submanifold (this roughly corresponds to equivariant coherent sheaves under mirror symmetry). The Floer cohomology of equivariant Lagrangian submanifolds has a natural endomorphism, which induces an $\R$-grading by generalized eigenspaces. Taking Euler characteristics with respect to the induced grading yields a deformation of the intersection number. Dehn twists act naturally on equivariant Lagrangians. Cotangent bundles and Lefschetz fibrations give fully computable examples. A key step in computations is to impose the ``dilation" condition stipulating that the BV operator applied to the symplectic cohomology class gives the identity.
\end{abstract}

\maketitle

\tableofcontents

\section{Introduction}
Let's begin with a somewhat schematic picture of Lagrangian Floer cohomology. Consider compact Lagrangian submanifolds $L$ in a symplectic manifold $M$ of dimension $2n$. Assuming that suitable assumptions have been imposed, one has a well-defined Floer cohomology group \cite{floer88c} $\HF^*(L_0,L_1)$ for any two such submanifolds $(L_0,L_1)$. For simplicity, we will work in the situation closest to classical cohomology, where Floer cohomology groups are $\Z$-graded and have coefficients in $\C$. They ``categorify'' the ordinary topological intersection number $L_0 \cdot L_1$, meaning that the Euler characteristic of Floer cohomology is
\begin{equation} \label{eq:int}
\chi(\HF^*(L_0,L_1)) = (-1)^{n(n+1)/2} L_0 \cdot L_1.
\end{equation}
As an example of how classical facts lifts to this level, let $\tau_V$ be the Dehn twist along a Lagrangian sphere $V\subset M$, (such symplectomorphisms arise naturally as Picard-Lefschetz monodromy maps). Then there is a long exact sequence \cite{seidel01}
\begin{equation}\label{eq:ps}
\xymatrix{
\HF^*(\tau_V(L_0),L_1) \ar[rr] && \HF^*(L_0,L_1) \ar[ld] \\
& \ar[ul]_{[1]} \mathit{Hom}(\HF^*(V,L_0), \HF^*(V,L_1)).
}
\end{equation}
Passing to the Euler characteristics recovers a version of the classical Picard-Lefschetz formula for the action of $\tau_V$ on homology.
``Categorification'' also brings in additional structures, not present at the intersection number level. Given any three Lagrangian submanifolds $(V,L_0,L_1),$ one has a composition map
\begin{equation}
\HF^*(L_1,L_2) \otimes \HF^*(L_0,L_1) \rightarrow \HF^*(L_0,L_2)
\end{equation}
(and indeed the $\swarrow$ in \eqref{eq:ps} is the dual of such a map). Moreover, any symplectomorphism $\phi : M \rightarrow M$ induces an isomorphism
\begin{equation} \label{eq:phi-star}
\phi_* : \HF^*(L_0,L_1) \rightarrow \HF^*(\phi(L_0),\phi(L_1)),
\end{equation}
which respects composition maps.

The aim of the present paper is to show that on certain noncompact symplectic manifolds, Lagrangian Floer cohomology theory can be further enriched by choosing an additional structure, called a {\em dilation}. This is a degree one element of the symplectic cohomology $\SH^*(M)$ in the sense of \cite{viterbo97a, cieliebak-floer-hofer95}, whose image under the BV (Batalin-Vilkovisky or loop rotation) operator is the canonical identity element in $\SH^0(M)$.

Fix $[b] \in \SH^1(M),$ not necessarily a dilation. We define the notion of a $b$-equivariant Lagrangian submanifold, which amounts to a Lagrangian submanifold satisfying the usual technical requirements of Floer cohomology together with a certain Lagrangian Floer cochain (see Definition~\ref{df:beq} for details). Any Lagrangian submanifold with vanishing first Betti number can be made $b$-equivariant, although this is not a necessary condition. Given two $b$-equivariant Lagrangian submanifolds, Section~\ref{sec:basic} shows that $\HF^*(L_0,L_1)$ can be equipped with a linear endomorphism $\tilde\Phi^1_{L_0,L_1}.$ Decomposing each graded piece of $\HF^*(L_0,L_1)$ into the generalized eigenspaces of $\tilde\Phi^1_{L_0,L_1}$ gives rise to a $\C$-grading on $\HF^*(L_0,L_1)$, which we call the equivariant grading. Floer continuation maps and composition maps respect the equivariant grading. Taking the Euler characteristic of each equivariant graded piece of $\HF^*(L_0,L_1)$ with respect to the usual grading on Floer cohomology, we obtain a function of a formal variable $q$ called the \emph{$q$-intersection number:}
\begin{equation} \label{eq:log-formula}
L_0 \bullet_q L_1 = \mathrm{Str}(e^{\log(q) \tilde{\Phi}^1_{L_0,L_1}}).
\end{equation}
Here, $\mathrm{Str}$ is the supertrace (Lefschetz trace). The Euler characteristic then turns out to be the $q = 1$ specialization of the $q$-intersection number, up to the same dimension-dependent sign as in \eqref{eq:int}. $q$-intersection numbers give a lower bound on the actual number of (transverse) intersection points of our two Lagrangian submanifolds, which sits in between the bound given by the ordinary intersection number and that provided by Floer cohomology itself. If $L$ is $b$-equivariant, Section~\ref{sec:dehn} shows that there is an induced equivariant structure on $\tau_V(L)$ such that the long exact sequence~\eqref{eq:ps} respects the equivariant grading. If $(L_0,L_1)$ are $b$-equivariant and $(\tau_V(L_0),\tau_V(L_1))$ are equipped with the induced equivariant structures, the isomorphism $\phi_*$ from \eqref{eq:phi-star} for $\phi = \tau_V$ respects the equivariant grading.

Even though this would violate some of the technical conditions in the body of the paper, one could in principle consider the same strategy for compact $M$, where symplectic cohomology reduces to ordinary cohomology, $\SH^*(M) \iso H^*(M;\C)$. However, then the whole reduces to classical topology: think of $b$ as represented by a $1$-form on $M$. A $b$-equivariant Lagrangian submanifold is a Lagrangian submanifold $L$ together with a function $c_L : L \rightarrow \C$ such that $dc_L = b|L$. Given $b$-equivariant Lagrangians $(L_0,c_{L_0}),(L_1,c_{L_1}),$ and $p \in L_0 \cap L_1,$ the grading of $[p] \in \HF^*(L_0,L_1)$ is $c_{L_0}(p) - c_{L_1}(p).$  The same holds for non-compact $M$, as long as the class $b$ belongs to the image of the canonical map $H^*(M) \rightarrow \SH^*(M).$

On the other hand, for general $b \in SH^1(M)$, $b$-equivariant Lagrangians do not appear to have an elementary interpretation. (In some examples, $q$-intersection numbers do yield the same result as previously known topological constructions \cite{givental88}.) The dilation condition which we impose singles out a special class of such $b$, which are guaranteed not to come from $H^*(M)$ (the BV operator kills classes coming from $H^*(M)$, but maps dilations to the identity, by definition). More importantly, this condition allows some computations of $q$-intersection numbers to be carried out just using the Poincar\'e duality properties of Floer cohomology. In particular, the equivariant grading of $\mathit{HF}^n(L,L) \iso \C$ is always $1.$ As a consequence of this and the long exact sequence \eqref{eq:ps}, we obtain the following \emph{$q$-Picard Lefschetz} formula:
\begin{equation}
\tau_{V}(L_0) \bullet_q L_1 = L_0 \bullet_q L_1 + (-1)^{n+1} q^{-1} (L_0 \bullet_q V) (V \bullet_q L_1).
\end{equation}

Obviously, the applicability of this theory depends on finding a dilation. When $M$ is a cotangent bundle, the relation with free loop space homology \cite{viterbo97b,salamon-weber03, abbondandolo-schwarz06} can be used to study this issue. See Section~\ref{sec:viterbo}. For instance, $M = T^*\C P^{n/2}$ carries a dilation. If we then take $L = \C P^{n/2}$ to be the zero-section, the action of $\tilde{\Phi}^1$ on $\HF^d(L,L) \iso H^d(L;\C)$ is $d/n$ times the identity, so that
\begin{equation} \label{eq:cotangent-pm}
L \bullet_q L = 1 + q^{2/n} + q^{4/n} + \cdots + q.
\end{equation}
Beyond the cotangent bundle case, one can use Lefschetz fibration methods, or other computational tools developed in the framework of symplectic cohomology theory, to prove that dilations exist in some cases. Section~\ref{sec:lefschetz} uses Lefschetz fibrations to show that the Milnor fibre of the $n$-dimensional $(A_m)$ singularity admits a dilation for $n \geq 3$ and $m \geq 1.$ However, a general characterization of manifolds admitting dilations is not available at present.

The original motivation for our construction, following suggestions of Bezrukavnikov and Okounkov, comes from Homological Mirror Symmetry \cite{kontsevich94}. Let $X$ be a quasi-projective variety with trivial canonical bundle, and $\xi$ an algebraic vector field on it. Let $E$ be an object of the bounded derived category $D^b_{\mathit{cpt}}(X)$ of compactly supported coherent sheaves on $X$. Then, $\xi$ gives rise to an infinitesimal deformation of $E$, which means a class in $\mathit{Ext}^1(E,E)$. We call $E$ infinitesimally equivariant if this class vanishes (and actually want to choose a coboundary for the relevant cocycle in an appropriate chain complex). If $E_0$ and $E_1$ are both infinitesimally equivariant, the action of $\xi$ defines a linear endomorphism of $\mathit{Ext}^*(E_0,E_1)$, which can be used to define a $q$-deformed version of the Mukai pairing. To understand how this is related to symplectic cohomology, observe that our vector field can be viewed as an element of Hochschild cohomology
\begin{equation}\label{eq:hkr}
(0,Z) \in H^1(X,\mathcal{O}_X) \oplus \Gamma(X,TX) \iso \HH^1(X).
\end{equation}
By definition \cite{kontsevich94}, the Hochschild cohomology $\HH^*(X)$ is $\Ext^*_{X \times X}(\mathcal{O}_\Delta,\mathcal{O}_\Delta)$ where $\mathcal{O}_\Delta$ denotes the structures sheaf of the diagonal $\Delta \subset X \times X.$ The isomorphism \eqref{eq:hkr} and its higher-degree versions generalize the Hochschild-Kostant-Rosenberg theorem, as pointed out in \cite{kontsevich94}. There is a canonical map from $\HH^*(X)$ to the Hochschild cohomology of the dg category underlying $D^b_{\mathit{cpt}}(X)$ (because of the support condition, we can't hope for this map to be an isomorphism, but see \cite{toen07} for an isomorphism result in a better-behaved framework). Now suppose that $M$ is the mirror of $X$. There is a natural map from $\SH^*(M)$ to the Hochschild cohomology of the Fukaya $A_\infty$-category \cite{seidel02,abouzaid10}. Hence, it seems natural to consider elements of $\SH^1(M)$ as formal analogues of algebraic vector fields. The dilation condition corresponds to asking for our vector field to contract the holomorphic volume form on $X$, at least in those cases where the mirror construction is based on such a form (as one can see from \cite{auroux} for instance, mirrors of open manifolds can also be constructed using volume forms with poles).

{\em Acknowledgments.} Roman Bezrukavnikov is the unofficial third author of this paper, and we are deeply indebted to his ideas. The first author would like to thank MSRI for its hospitality, and the NSF for partial financial support through grant DMS-0652620. The second author would like to thank the NSF for partial financial support through grant DMS-0703722 and the ISF for partial financial support through grant 1321/09. The second author was also partially supported by Marie Curie grant No. 239381. 

\section{The formal setup\label{sec:tqft}}

Our main arguments follow a version of the familiar TQFT framework, which involves a class of Riemann surfaces together with a fixed target space. We'll begin by outlining that framework in a semi-realistic language, which omits many necessary restrictions and details. The next section will contain a more precise description of the actual Floer-theoretic implementation.

Our target space will be a symplectic manifold $M$ of dimension $2n$, equipped with a Hamiltonian function $H$. For each $\lambda \in \R \setminus \{0\}$, we have a ``closed string'' chain complex $\CF^*(\lambda H)$ (all chain complexes are $\Z$-graded complexes of vector spaces over a fixed coefficient field $\K$, with differentials of degree $1$), whose differential is denoted by $d$. We consider Lagrangian submanifolds $L \subset M$. For every pair of such submanifolds, we have an ``open string'' chain complex $\CF^*(L_0,L_1)$, with differential $\mu^1 = \mu^1_{L_0,L_1}$.

The Riemann surfaces under consideration are of the form $S = \bar{S} \setminus \Sigma$, where $\bar{S}$ is a compact connected Riemann surface with (possibly) boundary, and $\Sigma \subset \bar{S}$ is a finite set of points. We exclude the case of a closed $S$. Each connected component $C \subset \partial S$ should come labeled with a Lagrangian submanifold $L_C$. Each $\zeta \in \Sigma \setminus \partial\bar{S}$ should be equipped with a distinguished tangent direction, by which we mean a half-line $\tau_\zeta \subset T_\zeta\bar{S}$. Additionally, the surface should come with a real one-form $\gamma \in \Omega^1(S)$, which satisfies
\begin{equation} \label{eq:sub-closed}
d\gamma \leq 0
\end{equation}
everywhere (with respect to the complex orientation), which vanishes near $\Sigma \cap \partial\bar{S}$, and which is closed near $\Sigma \setminus \partial\bar{S}$. We divide the marked points into inputs and outputs, $\Sigma = \Sigma_{\mathrm{in}} \cup \Sigma_{\mathrm{out}}$. Given $\zeta \in \Sigma \cap \partial\bar{S}$, let $(L_{\zeta,0},L_{\zeta,1})$ be the pair of Lagrangian submanifolds associated to the adjacent boundary components. More precisely, if $\zeta$ is an input, then $L_{\zeta,0}$ appears before $L_{\zeta,1}$ with respect to the natural boundary orientation; for outputs, the opposite holds. Similarly, given $\zeta \in \Sigma \setminus \partial\bar{S}$, let $\lambda_\zeta$ be the result of integrating $\gamma$ along a small circle around $\zeta$. If $\zeta$ is an output, this circle is oriented anticlockwise, and otherwise clockwise. We assume that $\lambda_\zeta \neq 0$ for all $\zeta$, so that $\CF^*(\lambda_\zeta H)$ is defined.

\begin{remark}
Our ``inputs'' correspond to ``outgoing ends'' in \cite{seidel04}, and our ``outputs'' to ``incoming ends''. This reversal of terminology brings it more in line with the rest of the literature, and with the algebraic function of the ends in Floer {\em co}homology theory.
\end{remark}

The simplest object, associated to a single surface $S$ as just described, is a chain map
\begin{equation} \label{eq:chain-map}
\xymatrix{\displaystyle
\qquad \bigotimes_{\zeta \in \Sigma_{\mathrm{in}} \setminus \partial \bar{S}} \CF^*(\lambda_\zeta H) \;\;\; \otimes \bigotimes_{\zeta \in \Sigma_{\mathrm{in}} \cap \partial \bar{S}} \CF^*(L_{\zeta,0},L_{\zeta,1}) \ar[d]^-{\phi_S} \\ \qquad \displaystyle
\Big(\bigotimes_{\zeta \in \Sigma_{\mathrm{out}} \setminus \partial \bar{S}} \CF^*(\lambda_\zeta H) \;\;\; \otimes \bigotimes_{\zeta \in \Sigma_{\mathrm{out}} \cap \partial \bar{S}} \CF^*(L_{\zeta,0},L_{\zeta,1})\Big)[\chi],
}
\end{equation}
where $\chi = - n \chi(\bar{S}) + 2n \,|\Sigma_{\mathrm{out}} \setminus \partial \bar{S}| + n\, |\Sigma_{\mathrm{out}} \cap \partial \bar{S}|$. Special cases are those of an infinite cylinder $\R \times S^1$ or strip $\R \times [0,1]$, with one input and one output, and equipped with a translation-invariant form $\gamma = \lambda \mathit{dt}$. In the first case, we require that the two tangent directions at the marked points of $\bar{S} = \C P^1$ should point towards each other. In the second case, we equip the surface with arbitrary boundary conditions $(L_0,L_1)$. The unitality axiom then requires that \eqref{eq:chain-map} should be equal to the identity on $\CF^*(\lambda H)$ and $\CF^*(L_0,L_1)$, respectively.

Crucially, the formalism extends to include families of surfaces $S$ parametrized by oriented manifolds. In the simplest case, if the parameter space is closed of some dimension $l$, the family gives rise to a chain map whose degree differs from the one given above by $-l$. More generally, we want to allow the parameter space to be a manifold with corners, and where the surfaces $S$ may degenerate along the boundary strata, by stretching them along tubular or strip-like ends. Such a family gives rise to a chain homotopy between the operations associated to its boundary faces. Instead of giving the abstract formulation, we prefer to just discuss the examples that are important for our purpose.

\begin{figure}
\begin{centering}
\includegraphics{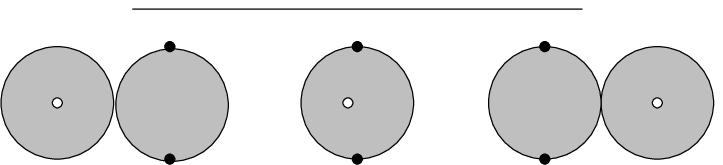}
\caption{\label{fig:h}}
\end{centering}
\end{figure}

We start with operations involving only the open string sector, which means no interior marked points. Here, the choice of one-form is essentially irrelevant, so we can just set it equal to zero. Take $\bar{S}$ to be a disc, with a single marked boundary point, and with $\partial S$ labeled by some $L$. By considering the marked point either as an output or input, we get
\begin{equation} \label{eq:relative-unit}
\begin{aligned}
& e_L \in \CF^0(L,L), \\
& e_L^\vee: \CF^n(L,L) \longrightarrow \K.
\end{aligned}
\end{equation}
Next, take the same $\bar{S}$ with three marked boundary points, of which two are inputs and one is an output. If we label the three boundary components of the resulting $S$ with Lagrangian submanifolds $L_0,L_1,L_2$, we get a chain map
\begin{equation} \label{eq:multiplication}
\mu^2 = \mu^2_{L_0,L_1,L_2}: \CF^*(L_1,L_2) \otimes \CF^*(L_0,L_1) \longrightarrow \CF^*(L_0,L_2).
\end{equation}
By considering the universal one-parameter family of discs with four marked boundary points, one obtains chain homotopies which show that $\mu^2$ is homotopy associative. For the sake of familiarity, we keep the sign conventions close to those for ordinary differential graded algebras (see Section \ref{sec:signs} for further discussion); so the homotopy has the form
\begin{equation} \label{eq:mu3}
\begin{aligned}
& \mu^3 = \mu^3_{L_0,L_1,L_2,L_3}: \CF^*(L_2,L_3) \otimes \CF^*(L_1,L_2) \otimes \CF^*(L_0,L_1) \\ & \qquad \qquad \qquad \qquad \longrightarrow \CF^*(L_0,L_3)[-1], \\
& \mu^1(\mu^3(a_3,a_2,a_1)) + \mu^3(\mu^1(a_3),a_2,a_1) + (-1)^{|a_3|} \mu^3(a_3,\mu^1(a_2),a_1) \\ & +
(-1)^{|a_3|+|a_2|} \mu^3(a_3,a_2,\mu^1(a_1)) = \mu^2(a_3,\mu^2(a_2,a_1)) - \mu^2(\mu^2(a_3,a_2),a_1).
\end{aligned}
\end{equation}
This is well-known, being the starting point for the construction of Fukaya $A_\infty$-structures. We will not pursue this further, and instead turn to the relation with \eqref{eq:relative-unit}. Using a one-parameter family of surfaces, one shows that $e_L$ is a two-sided unit up to homotopy for $\mu^2$. A related weaker statement, which does not require the unitality axiom, says that $e_L$ is a central element in cohomology. The corresponding statement for $e_L^\vee$ says that the pairing obtained by composing it with $\mu^2$, namely
\begin{equation} \label{eq:pairing}
\begin{aligned}
& \CF^*(L_1,L_0) \otimes \CF^*(L_0,L_1) \longrightarrow \K[-n], \\
& (a_2,a_1) \longmapsto e_{L_0}^\vee(\mu^2(a_2,a_1))
\end{aligned}
\end{equation}
is (graded) symmetric up to homotopy. This time, it is worth while spelling out the details, since the homotopy involved will appear in our computations later on. It is of the form
\begin{equation} \label{eq:hvee}
\begin{aligned}
& h_{L_0,L_1}^\vee: \CF^*(L_1,L_0) \otimes \CF^*(L_0,L_1) \longrightarrow \K[-n-1], \\
& h_{L_0,L_1}^\vee(\mu^1(a_2),a_1) + (-1)^{|a_2|} h_{L_0,L_1}^\vee(a_2,\mu^1(a_1)) \\ & = e_{L_0}^\vee(\mu^2(a_2,a_1)) - (-1)^{|a_1|\,|a_2|} e_{L_1}^\vee(\mu^2(a_1,a_2))
\end{aligned}
\end{equation}
and arises from a one-parameter family of discs with two marked boundary points, both of which are considered as inputs. To see that family in more familiar terms, it is convenient to introduce an additional interior marked point; this point is not an input or output, it just serves as a marker which partially breaks the symmetry of the disc. So, we consider the moduli space of discs with one interior and two boundary marked points, and its Deligne-Mumford compactification. This is a copy of $\bar{\R} = \R \cup \{\pm\infty\}$, where the points at infinity correspond to surfaces with two components (see Figure \ref{fig:h}, where the additional marked point is drawn as a white dot). Those surfaces indeed correspond to the two compositions which appear on the right side of the equation in \eqref{eq:hvee}.

We next turn to operations involving only the closed string sector, which means surfaces with no boundary. The simplest case, already mentioned above, is $S = \R \times S^1$ with $\gamma = \lambda \mathit{dt}$, with one input and one output. As a single surface, this just yields the identity map on $\CF^*(\lambda H)$, but one can get a family of surfaces parametrized by $S^1$ by rotating the tangent direction at one of the marked points (because of the automorphism group, it is irrelevant which one we pick). This gives rise to a degree $-1$ chain map, the Batalin-Vilkovisky operator
\begin{equation} \label{eq:bv}
\delta: \CF^*(\lambda H) \longrightarrow \CF^*(\lambda H)[-1].
\end{equation}
One can show by using appropriate two-parameter families that $\delta^2$ is nullhomotopic.

Take $\bar{S}$ to be a sphere with a single point marked as output, but allowing for non-closed $\gamma$. This is subject to the asymmetry condition \eqref{eq:sub-closed}, so the resulting element
\begin{equation} \label{eq:closed-unit}
\begin{aligned}
& e \in \CF^0(\lambda H)
\end{aligned}
\end{equation}
only exists for $\lambda > 0$. Of course $d(e) = 0$, and beyond that, a suitable two-parameter family argument shows that $\delta e$ is a $d$-coboundary. Next, we return to the case of a sphere with one input and one output, first fixing the tangent directions at both points, but allowing a general $\gamma$. This gives rise to continuation maps
\begin{equation} \label{eq:continuation}
\CF^*(\lambda_0 H) \longrightarrow \CF^*(\lambda_1 H)
\end{equation}
for all $\lambda_0 \leq \lambda_1$ (in the case $\lambda_0 = \lambda_1$, this is again the identity). One can show that up to chain homotopy, these compose well and are compatible with \eqref{eq:bv} as well as \eqref{eq:closed-unit}. A standard use of these maps is to remove the dependence on the parameter $\lambda$ by passing to direct or inverse limits, see Section \ref{sec:viterbo} below. This idea can also be combined with the one underlying the definition of $\delta$.

Open-closed string maps are obtained by looking at surfaces which involve both kinds of marked points. The simplest example is to take $\bar{S}$ to be a disc with one interior point and one boundary point. The interior marked point is an input, and the tangent direction and that point goes towards the boundary marked point. The boundary marked point can be either an output or an input, leading to chain maps
\begin{equation} \label{eq:phi0}
\begin{aligned}
& \phi^0_L: \CF^*(\lambda H) \longrightarrow \CF^*(L,L), \\
& \phi^{0,\vee}_L: \CF^*(\lambda H) \otimes \CF^{n-*}(L,L) \longrightarrow \K,
\end{aligned}
\end{equation}
both defined for all $\lambda \neq 0$. Suppose that $\lambda > 0$, and consider $e \in \CF^0(\lambda H)$. A suitable one-parameter family argument shows that $\phi^0_L(e) - e_L$ as well as $\phi^{0,\vee}_L(e,\cdot) - e_L^\vee$ are coboundaries.

\begin{figure}
\begin{centering}
\includegraphics{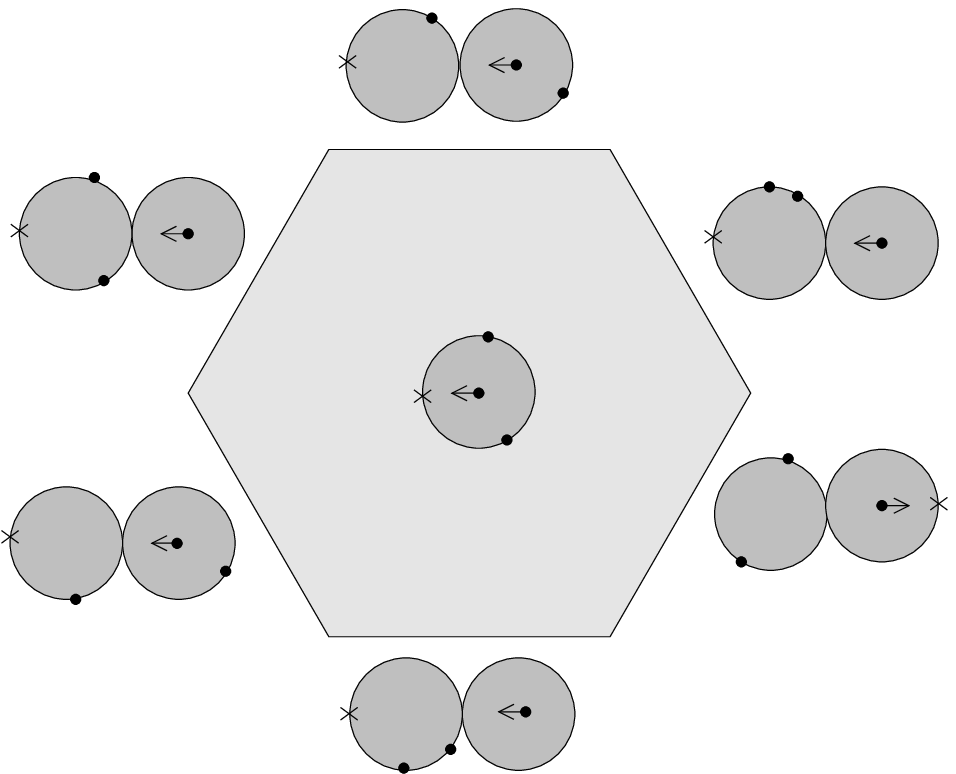}
\caption{\label{fig:hex}}
\end{centering}
\end{figure}

Let's still take $\bar{S}$ to be a disc, but now with one interior input point and two boundary marked points, one an input and one an output, and where the tangent line points towards the output boundary point. This leads to essentially the same one-parameter family as in Figure \ref{fig:h}, which yields a homotopy
\begin{equation} \label{eq:phi1}
\begin{aligned}
& \phi^1_{L_0,L_1}: \CF^*(\lambda H) \otimes \CF^*(L_0,L_1) \longrightarrow \CF^*(L_0,L_1)[-1], \\
& \mu^1(\phi^1_{L_0,L_1}(b,a)) + \phi^1_{L_0,L_1}(db,a) + (-1)^{|b|} \phi^1_{L_0,L_1}(b,\mu^1(a)) \\ & = \mu^2(\phi^0_{L_1}(b),a) - (-1)^{|a|\,|b|} \mu^2(a,\phi^0_{L_0}(b)).
\end{aligned}
\end{equation}
In the next step, one would take one interior and three boundary marked points, where again the tangent direction points towards the unique boundary marked point which is an output. The corresponding compactified moduli space (shown in Figure \ref{fig:hex}; the point marked with a cross is the output) is a hexagon. It leads to a map
\begin{equation} \label{eq:phi2}
\begin{aligned}
& \phi^2_{L_0,L_1,L_2}: \CF^*(\lambda H) \otimes \CF^*(L_1,L_2) \otimes \CF^*(L_0,L_1) \longrightarrow \CF^*(L_0,L_2)[-2], \\
& \mu^1(\phi^2_{L_0,L_1,L_2}(b,a_2,a_1)) - \phi^2_{L_0,L_1,L_2}(db,a_2,a_1) \\ & - (-1)^{|b|} \phi^2_{L_0,L_1,L_2}(b,\mu^1(a_2),a_1) - (-1)^{|b|+|a_2|} \phi^2_{L_0,L_1,L_2}(b,a_2,\mu^1(a_1))\\ & =
-\mu^3(\phi^0_{L_2}(b),a_2,a_1) + (-1)^{|a_2|\,|b|} \mu^3(a_2,\phi^0_{L_1}(b),a_1) \\ & - (-1)^{(|a_2|+|a_1|)|b|} \mu^3(a_2,a_1,\phi^0_{L_0}(b))
+ \phi^1_{L_0,L_2}(b,\mu^2(a_2,a_1)) \\ & -
\mu^2(\phi^1_{L_1,L_2}(b,a_2),a_1) - (-1)^{(|b|+1)|a_2|} \mu^2(a_2,\phi^1_{L_0,L_1}(b,a_1)).
\end{aligned}
\end{equation}
Again, these are just the first of an infinite sequence of maps $\phi^k_{L_0,\dots,L_k}$, which have a natural interpretation in terms of the Hochschild cochain complex of the Fukaya $A_\infty$-category, see \cite{seidel02} or the more recent \cite{abouzaid10}.

\begin{figure}
\begin{centering}
\includegraphics{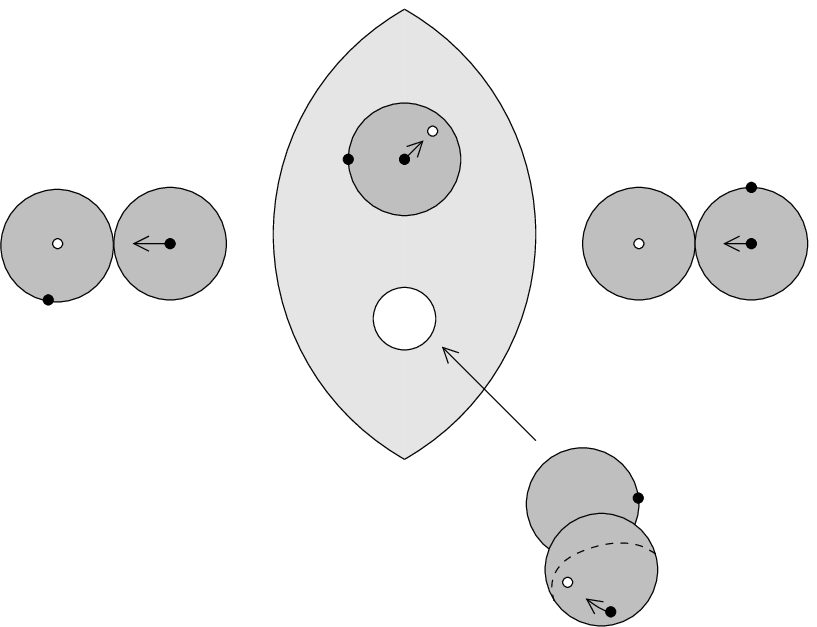}
\caption{\label{fig:k}}
\end{centering}
\end{figure}

The last piece of structure which is of direct importance for our purpose measures the interaction of $\phi^1_{L,L}$ and $e^\vee_L$, and involves the BV operator $\delta$. It is of the form
\begin{equation} \label{eq:kvee}
\begin{aligned}
& k_L^\vee : \CF^*(\lambda H) \otimes \CF^*(L,L) \longrightarrow \K[-n-2], \\
& k_L^\vee(db,a) + (-1)^{|b|} k_L^\vee(b,\mu^1(a)) \\ & =
e_L^\vee(\phi^1_{L,L}(b,a)) - \phi^{0,\vee}_L(\delta b,a) - h_{L,L}^\vee(\phi^0_L(b),a).
\end{aligned}
\end{equation}
To define this, we consider the moduli space of discs with two marked interior points and one marked boundary point. One of the interior points is an input, and the tangent line points towards the other interior point, which just serves as a marker. The boundary point is also an input. The moduli space itself is an open annulus; the compactification which plays a role here is the ``real blowup'' of the Deligne-Mumford compactification from \cite{kimura-stasheff-voronov95}, shown in Figure \ref{fig:k}.

In that space, the ``blown up'' boundary circle associated to bubbling off of a sphere leads to the term in \eqref{eq:kvee} involving $\delta$. More precisely, it parametrizes nodal surfaces
\begin{equation}
\bar{S} = \bar{S}_+ \cup \bar{S}_-
\end{equation}
where $\bar{S}_+$ is a disc with one interior marked point $\zeta_{+,0}$ and one boundary marked point $\zeta_{+,1}$, and $\bar{S}_-$ is a sphere with three marked points $\zeta_{-,0}$, $\zeta_{-,1}$, $\zeta_{-,2}$ (the first is the node, and the last is the ``marker''). This comes with a distinguished tangent direction at $\zeta_{-,1}$ pointing towards $\zeta_{-,2}$, and in addition, by definition of the compactification, with an ``gluing angle'', which is a real half-line inside
\begin{equation} \label{eq:gluing-angle}
(T\bar{S}_+)_{\zeta_{+,0}} \otimes_\C (T\bar{S}_-)_{\zeta_{-,0}}.
\end{equation}
One can identify $\bar{S}_+$ with the standard disc in a preferred way, so that the marked points go to $0$ and $1$, respectively. Having done that, there is a canonical tangent direction at $\zeta_{+,0}$, namely the one pointing towards $\zeta_{+,1}$. This and \eqref{eq:gluing-angle} together give rise to a distinguished tangent direction at $\zeta_{-,0}$. One can identify $\bar{S}_- \iso \C P^1$ in such a way that the marked points go to $0,\infty$ and $1$, respectively. It then carries distinguished tangent directions at $0$ (variable) and $\infty$ (fixed, pointing towards $1$), hence can be identified with the family of surfaces which occur in the definition of $\delta$.

So far, we have rigidly upheld the distinction between inputs and outputs, but in the Floer-theoretic framework this will turn out to be somewhat artificial, if technically convenient. Let's introduce the standard notation for the Floer cohomology groups arising from our complexes, $\HF^*(\lambda H) = H^*(\CF(\lambda H))$ and $\HF^*(L_0,L_1) = H^*(\CF(L_0,L_1))$. The pairing \eqref{eq:pairing} is nondegenerate on the cohomology level, hence induces a Poincar{\'e} duality isomorphism
\begin{equation} \label{eq:poincare}
\HF^i(L_1,L_0) \iso \HF^{n-i}(L_0,L_1)^\vee,
\end{equation}
which allows one to turn inputs into outputs, and vice versa. For instance, for $L_0 = L_1 = L$, it maps $[e_L]$ to $[e_L^\vee]$, and also relates the cohomology level maps induced by $\phi^0_L$ and $\phi^{0,\vee}_L$. There is a similar duality for closed strings, $\HF^*(-\lambda H) \iso \HF^{2n-*}(\lambda H)^\vee$, but that will play no role in our applications.

\section{The geometric setup\label{sec:analysis}}

The kind of framework sketched out above can be implemented in Floer theory in several related but different ways. For our intended applications, it is crucial that the BV operator be nontrivial, which means that our symplectic manifolds should be noncompact. Other than that, we remain in the technically simplest situation, where strong exactness and Calabi-Yau type conditions are imposed.

\begin{assumption} \label{th:convexity}
Let $M^{2n}$ be a symplectic manifold which is exact, $\omega = d\theta$, together with a compatible almost complex structure $I$. We assume that $c_1(M) = 0$, and in fact choose a trivialization of the canonical bundle $\scrK_M = \Lambda^n_\C(T^*M)$ (only the homotopy class of the trivialization will be important). In addition, we choose a class $[\alpha] \in H^2(M;\Z/2)$, and a $\Z/2$-gerbe $\alpha$ representing that class. 

The other piece of data we need to fix is a Hamiltonian $H \in \smooth(M,\R)$. Let $X$ be its Hamiltonian vector field. We also need to assume some control over its geometry at infinity. Namely, $M$ should admit an exhaustion (an increasing sequence of relatively compact open subsets $U_1 \subset U_2 \subset \cdots$, whose union is $M$) such that the following holds. Suppose that $S$ is a connected compact Riemann surface with nonempty boundary, equipped with a two-form $\gamma$ satisfying $d\gamma \leq 0$, and $u: S \rightarrow M$ a solution of
\begin{equation} \label{eq:standard-equation}
(du - X \otimes \gamma)^{0,1} = \half (du + I \circ du \circ i - X \otimes \gamma - IX \otimes \gamma \circ i) = 0
\end{equation}
such that $u(\partial S) \subset U_k$ for some $k$. Then $u(S) \subset U_{k+1}$.
\end{assumption}

Let ${\mathcal P} \subset \R$ be the set of all $\lambda$ such that the one-periodic orbits of $\lambda X$ are not contained in any compact subset of $M$. This is the set of forbidden values: for $\lambda \in {\mathcal P}$, $\CF^*(\lambda H)$ is not defined (this is the first detail in which we diverge from the formal setup explained in Section \ref{sec:tqft}). For all other values, one proceeds as follows. Choose a time-dependent $K_d \in \smooth(S^1 \times M,\R)$ satisfying $K_d(t,x) = \lambda H(x)$ for all $x$ lying outside some compact subset of $M$ (the subscript $d$ refers to the Floer differential). This should be such that the one-periodic orbits of the associated time-dependent vector field $X_d$ are nondegenerate. By definition of ${\mathcal P}$, there are necessarily only finitely many of them. $\CF^*(\lambda H)$ has one generator for each such orbit $y = y(t)$, with degrees given by Conley-Zehnder indices (our normalization is such that for small time-independent functions, the Conley-Zehnder index agrees with the Morse index). Similarly, one takes a family $J_d$ of almost complex structures on $M$ depending on $t \in S^1$, such that $J_d(t,x) = I(x)$ outside a compact subset. The differential $d$ on $\CF^*(\lambda H)$ is obtained by counting solutions of the associated Floer equation
\begin{equation} \label{eq:floer}
\left\{ \begin{aligned} & u: \R \times S^1 \longrightarrow M, \\
& \partial_s u + J_d(t,u)(\partial_t u - X_d(t,u)) = 0, \\
& \textstyle\lim_{s \rightarrow \pm \infty} u(s,t) = y_\pm(t).
\end{aligned}
\right.
\end{equation}
We remind the reader that this is formally a negative gradient flow equation for the action functional
\begin{equation} \label{eq:action}
A(y) = \textstyle\int_{S^1} -y^*\theta + K_d(t,y(t)) \mathit{dt}.
\end{equation}
The (cohomological) convention here is that solutions with limits $y_\pm$ contribute to the coefficient of $y_-$ in $d(y_+)$. Outside a compact subset of $M$, \eqref{eq:floer} is of the form \eqref{eq:standard-equation}, hence we know from Assumption \ref{th:convexity} that solutions cannot escape to infinity. The gerbe $\alpha$ determines a $\Z/2$-cover of the free loop space of $M$, and we use the associated $\K$-coefficient system to determine the signs with which solutions of \eqref{eq:floer} contribute to the differential (in the language of \cite{floer-hofer93}, different choices of $[\alpha]$ correspond to different coherent orientations). The same will hold for all the moduli spaces which occur later on. 

\begin{assumption} \label{th:lagrangian}
When considering Lagrangian submanifolds, we allow only those $L \subset M$ which are closed, connected, and exact (meaning that $\theta|L$ is an exact one-form). In addition, we assume that each $L$ comes equipped with a grading \cite{seidel99}, which in particular induces an orientation, and with a $\mathit{Spin}$ structure relative to $\alpha$ \cite{fooo}. 
\end{assumption}

The grading and relative $\mathit{Spin}$ structure allow us to have integrally graded Floer cohomology groups $\mathit{HF}^*(L_0,L_1)$ with coefficients in an arbitrary $\K$. To define that for a given pair $(L_0,L_1)$, choose a time-dependent $K_{\mu^1} \in \smooth([0,1] \times M,\R)$ which is compactly supported. This should be such that the length one chords of the associated time-dependent vector field $X_{\mu^1}$ going from $L_0$ to $L_1$ are nondegenerate. $\CF^*(L_0,L_1)$ has one generator for each such chord, with degrees given by the absolute Maslov indices, that in turn depend on the gradings. Similarly, one takes a family $J_{\mu^1}$ of almost complex structures on $M$ depending on $t \in [0,1]$, such that $J_{\mu^1}(t,x) = I(x)$ outside a compact subset. The differential $\mu^1$ on $\CF^*(L_0,L_1)$ is obtained by counting solutions of an equation like \eqref{eq:floer}, where now the domain is $\R \times [0,1]$, and we have boundary conditions $u(s,0) \in L_0$, $u(s,1) \in L_1$. In the special case $L_0 = L_1 = L$, we have the PSS isomorphism \cite{piunikhin-salamon-schwarz94,albers08}
\begin{equation} \label{eq:pss}
\mathit{HF}^*(L,L) \iso H^*(L;\K).
\end{equation}
Jumping ahead slightly, we should mention that with respect to this isomorphism, $[e_L]$ is the identity, $[e_L^\vee]$ the fundamental class, and $[\mu^2]$ is the ordinary cup-product.

We now quickly sketch the setup that defines operations associated to general Riemann surfaces. Let's start with the case of a single surface $S$, as in \eqref{eq:chain-map}, but with the additional assumption that $\lambda_\zeta \notin {\mathcal P}$ for all $\zeta \in \Sigma \setminus \partial \bar{S}$. Choose distinguished local holomorphic coordinates near each marked point $\zeta \in \Sigma$, of the form
\begin{equation}
\begin{cases}
\epsilon_\zeta: \R^+ \times [0,1] \rightarrow S, & \text{for $\zeta \in \Sigma_{\mathrm{in}} \cap \partial \bar{S}$,} \\
\epsilon_\zeta: \R^- \times [0,1] \rightarrow S, & \text{for $\zeta \in \Sigma_{\mathrm{out}} \cap \partial \bar{S}$,} \\
\epsilon_\zeta: \R^+ \times S^1 \rightarrow S, & \text{for $\zeta \in \Sigma_{\mathrm{in}} \setminus \partial \bar{S}$,} \\
\epsilon_\zeta: \R^- \times S^1 \rightarrow S, & \text{for $\zeta \in \Sigma_{\mathrm{out}} \setminus \partial \bar{S}$.}
\end{cases}
\end{equation}
In the last two cases, we assume that the preferred tangent direction at the marked point is that tangent to the half-line $\epsilon_\zeta(\R^\pm \times \{0\})$. We may assume that $\epsilon_\zeta^*\gamma = 0$ for boundary marked points $\zeta \in \Sigma \cap \partial \bar{S}$. A minor technical modification of the given $\gamma$ ensures that similarly, $\epsilon_\zeta^*\gamma = \lambda_\zeta \mathit{dt}$ for interior marked points $\gamma \in \Sigma \setminus \partial \bar{S}$, which we assume to be the case from now on.

\begin{definition}
A perturbation datum is a pair $(K,J)$ of the following form. $K$ is a one-form on $S$ with values in the space of functions on $M$, or equivalently a section of the pullback bundle $T^*S \rightarrow S \times M$. $J$ is a family of compatible almost complex structures on $M$ parametrized by $S$. They should have the following properties.

Outside a compact subset of $M$ we have $K = H \otimes \gamma$ and $J = I$. If $L_C$ is the Lagrangian submanifold associated to some component $C \subset \partial S$, and $\xi \in TC$ is a vector tangent to that component, then $K(\xi)|L_C = 0$. Over a strip-like end associated to an interior marked point we want to have $(K,J) = (K_d \,\mathit{dt},J_d)$, as in the definition of Floer cohomology with the constant $\lambda = \lambda_\zeta$. Similarly, over the ends associated to boundary marked points, we have $(K,J) = (K_{\mu^1}\, \mathit{dt},J_{\mu^1})$, as in the definition of Floer cohomology for the pair $(L_{\zeta,0},L_{\zeta,1})$.
\end{definition}

Given these data, one can consider solutions of
\begin{equation} \label{eq:eq}
\left\{
\begin{aligned}
& u: S \longrightarrow M, \\
& u(z) \in L_C \text{ for $z \in C \subset \partial S$}, \\
& (du - Y)^{0,1} = 0,
\end{aligned}
\right.
\end{equation}
where $Y$ is the section of $\Hom_\R(TS,TM) \rightarrow S \times M$ such that for each $\xi \in TS$, $Y(\xi)$ is the Hamiltonian vector field of $K(\xi)$. At every end, we require that $u(\epsilon_\zeta(s,t)) \rightarrow y_\zeta(t)$, where $y_\zeta$ corresponds to a generator of the relevant Floer cochain complex. Counting solutions of \eqref{eq:eq} with these asymptotics yields a chain map \eqref{eq:chain-map}.

The generalization to families of Riemann surfaces over a closed parameter space is relatively straightforward. Due to our requirement that the interior marked points should come with distinguished tangent lines, one can choose strip-like ends continuously over any such family, in a way which is unique up to homotopy. On then chooses a suitable family of data $(K,J)$ whose restriction to the strip-like ends is the same as before, hence constant in the family, and considers the parametrized version of \eqref{eq:eq}. The example which is most relevant for us is the BV operator \eqref{eq:bv}. This involves the constant family of cylinders $S = \R \times S^1$ but with strip-like ends which depend on a parameter $r \in S^1$, concretely $\epsilon_-(r,s,t) = (s,t)$ for $s \ll 0$, and $\epsilon_+(r,s,t) = (s,t+r)$ for $s \gg 0$. As before, we write $(K_d,J_d)$ for the data used to define $\CF^*(\lambda H)$. Choose a $K_{\delta} \in \smooth(S^1 \times \R \times S^1 \times M, \R)$ and a family $J_{\delta}$ of almost complex structures depending on $S^1 \times \R \times S^1$, such that
\begin{equation}
\left\{
\begin{aligned}
& K_{\delta}(r,s,t,x) = K_d(t,x), \;  J_{\delta}(r,s,t,x) = J_d(t,x) \text{ for $s \ll 0$}; \\
& K_{\delta}(r,s,t,x) = K_d(t+r,x), \; J_{\delta}(r,s,t,x) = J_d(t+r,x) \text{ for $s \gg 0$}; \\
& K_{\delta}(r,s,t,x) = \lambda H(x), \; J_{\delta}(r,s,t,x) = I(x) \text{ for $x$ outside a compact subset.}
\end{aligned}
\right.
\end{equation}
In terms of the previous general setup, this means that we look at the section $K_{\delta}\, \mathit{dt}$ of $\Hom_\R(TS,TM) \rightarrow S^1 \times S \times M$. The associated parametrized version of \eqref{eq:eq} looks concretely like this:
\begin{equation} \label{eq:bv-equation}
\left\{ \begin{aligned}
& r \in S^1, \\
& u: \R \times S^1 \longrightarrow M, \\
& \partial_s u + J_{\delta}(r,s,t,u)(\partial_t u - X_\delta(r,s,t,u)) = 0, \\
& \textstyle \lim_{s \rightarrow -\infty} u(s,t) = y_-(t), \\
& \textstyle \lim_{s \rightarrow +\infty} u(s,t) = y_+(t+r).
\end{aligned}
\right.
\end{equation}
In terms of the original construction of continuation maps \cite{salamon-zehnder92}, equations like \eqref{eq:bv-equation}, with the parameter valued in $r \in [0,1]$ rather than $S^1$, define chain homotopies between the continuation maps associated to the endpoints $r = 0,1$. In our case, the data at the endpoints agree, which is why we get a chain map of degree $-1$ instead.

\begin{remark} \label{th:equivariant}
In fact, there is an infinite sequence of operations $d_k$ of degree $1-2k$, with $d_0 = d$, $d_1 = \delta$. The higher order maps are also all defined in terms of families of continuation map equations. For instance, $d_2$ is the homotopy between $\delta^2$ and zero. The totality of higher order relations are most easily described by writing a series
\begin{equation} \label{eq:equivariant-differential}
d_{\eq} = d_0 + u d_1 + \cdots
\end{equation}
in a degree $2$ formal variable $u$, which then satisfies $d_{\eq}^2 = 0$. In fact, this is precisely the differential used to define $S^1$-equivariant Floer cohomology (see again \cite{viterbo97a} or \cite{seidel07}; for additional foundational material about equivariant Floer cohomology theories see \cite{hutchings08, seidel-smith10}; an analytically more demanding approach is presented in \cite{bourgeois-oancea09}). Take the ring $\K[[u]]$, and the $\K[[u]]$-module $\K[u^{-1}] = \K((u))/u\K[[u]]$. The equivariant chain complex is
\begin{equation}
\CF_{\eq}^*(\lambda H) = \CF^*(\lambda H) \otimes \K[u^{-1}],
\end{equation}
and it carries the $\K[[u]]$-linear differential defined by \eqref{eq:equivariant-differential}.
\end{remark}

In the case where the parameter space of the family is itself noncompact, one has to carefully control the limiting behaviour of $(K,J)$, so as to get moduli spaces with the appropriate behaviour at the boundary (this is referred to as ``consistency of choices'' in \cite{seidel04}). The example which is most relevant for us is \eqref{eq:kvee}, and in particular the appearance of the BV operator in it. First, $\phi_L^{0,\vee}$ can be defined by a choice of perturbation datum $(K_{\phi^\vee} \mathit{dt},J_{\phi^\vee})$ on the single surface $S = (\R^+ \times S^1) \setminus \{(0,0)\}$, which over the interior end equals $(K_d\, \mathit{dt},J_d)$. Define a family of such data, parametrized by $(q,r) \in [Q,\infty) \times S^1$ for some $Q \gg 0$, by gluing together this piece with the one defining $\delta$:
\begin{equation}
(K_{k^\vee}(q,r,s,t,x) \,\mathit{dt},J_{k^\vee}(q,r,s,t,x)) = \begin{cases} (K_{\phi^\vee}(s,t,x) \, \mathit{dt},J_{\phi^\vee}(s,t,x)) & s \leq q, \\
(K_\delta(r,s-2q,t,x) \,\mathit{dt},J_\delta(r,s-2q,t,x)) & s \geq q. \end{cases}
\end{equation}
This yields the perturbation data in a neighbourhood of the boundary circle in the moduli space from Figure \ref{fig:k}. More precisely, to identify our surface with the one depicted there, one would equip it with the additional marked point $(2q,r)$. By definition, the neck-stretching limit $q \rightarrow \infty$ decomposes our surface into two parts, which give rise to the expressions $\phi_L^{0,\vee}$ and $\delta$ in \eqref{eq:kvee}.

\section{The basic construction\label{sec:basic}}

We continue in the framework set out in Assumptions \ref{th:convexity} and \ref{th:lagrangian}. Fix some $\lambda \in \R^+ \setminus {\mathcal P}$, and let $e \in \CF^0(\lambda H)$ be the identity element.

\begin{definition}
A {\em dilation} is a cocycle $b \in \CF^1(\lambda H)$ satisfying
\begin{equation} \label{eq:liouville-element}
[\delta b] = [e] \in \HF^0(\lambda H).
\end{equation}
\end{definition}

We postpone the discussion of the existence of dilations to Section \ref{sec:viterbo}, and assume for now that such an element has been fixed.

\begin{definition}\label{df:beq}
A {\em $b$-equivariant Lagrangian submanifold} is an $L \subset M$ together with an equivalence class of cochains $c_L \in \CF^0(L,L)$ satisfying
\begin{equation}
\mu^1(c_L) = \phi^0_L(b).
\end{equation}
Here, two $c_L$ are considered equivalent if their difference is a degree zero coboundary.
\end{definition}

The obstruction to the existence of $c_L$ for a given $L$ is obviously $[\phi^0_L(b)] \in \HF^1(L,L) \iso H^1(L;\K)$. If that vanishes, the set of equivalence classes of choices is an affine space over $\HF^0(L,L) \iso H^0(L;\K) = \K$. We write
\begin{equation}
L \longmapsto L\langle s \rangle
\end{equation}
for the operation on $b$-equivariant Lagrangian submanifolds which subtracts $s$ times $e_L$ from a given $c_L$. This should not be confused with the more common shift operation $L \mapsto L[r]$, which changes the grading \cite[Section 11k]{seidel04}.

Suppose that $L_0$ and $L_1$ are $b$-equivariant. We can then define an endomorphism
\begin{equation} \label{eq:corrected-phi}
\begin{aligned}
& \tilde{\phi}^1_{L_0,L_1}: \CF^*(L_0,L_1) \longrightarrow \CF^*(L_0,L_1), \\
& \tilde{\phi}^1(a) = \phi^1(b,a) - \mu^2(c_{L_1},a) + \mu^2(a,c_{L_0}).
\end{aligned}
\end{equation}
It follows from \eqref{eq:phi1} that this is a chain map, hence induces an endomorphism of cohomology, which we denote by $\tilde\Phi^1_{L_0,L_1}$.

\begin{definition}
Suppose that our coefficient field is $\K = \C$. The endomorphism $\tilde\Phi^1_{L_0,L_1}$ determines a decomposition of $\HF^*(L_0,L_1)$ into (generalized) eigenspaces $\HF^*(L_0,L_1)^\lambda$, indexed by $\lambda \in \C$. The $q$-intersection number of two $b$-equivariant Lagrangian submanifolds $(L_0,L_1)$ is the function of one formal variable $q$ given by combining the Euler characteristic of the eigenspaces:
\begin{equation}
L_0 \bullet_q L_1 = \sum_\lambda \chi(\HF^*(L_0,L_1)^\lambda) q^\lambda.
\end{equation}
\end{definition}

This is easily seen to the expression \eqref{eq:log-formula} given in the Introduction. $\tilde{\Phi}^1_{L_0,L_1}$ depends only on the equivalence classes of $c_{L_0}$ and $c_{L_1}$. By definition it satisfies
\begin{equation} \label{eq:add-identity}
\tilde{\Phi}_{L_0\langle s_0 \rangle, L_1 \langle s_1 \rangle}^1 = \tilde{\Phi}_{L_0,L_1}^1 + (s_1 - s_0) \mathrm{id},
\end{equation}
which implies that
\begin{equation} \label{eq:eq-shift}
L_0\langle s_0 \rangle \bullet_q L_1\langle s_1 \rangle = q^{s_1-s_0} L_0 \bullet_q L_1.
\end{equation}

\begin{remark} \label{th:evaluation}
For $L_0 = L_1 = L$, $\tilde{\Phi}^1_{L,L}$ is entirely independent of the choice of $c_L$. In fact, using the PSS isomorphism \eqref{eq:pss}, one can get a more direct description of this element, which generalizes to not necessarily $b$-equivariant Lagrangian submanifolds. Namely, take the Riemann surface $S = \R^+ \times S^1$, with a perturbation datum which equals $(H_d \mathit{dt}, J_d)$ over the end. Consider the moduli space $\scrM(y)$ of solutions of the associated equation \eqref{eq:eq} with given limit $y$ and boundary conditions $u(\{0\} \times S^1) \subset L$. Take the double evaluation map on its compactification, given by
\begin{equation}
\begin{aligned}
& S^1 \times \bar\scrM(y) \longrightarrow L \times L, \\
& (r,u) \longmapsto (u(0,0),u(0,r)).
\end{aligned}
\end{equation}
The formal sum of such maps corresponding to the coefficients of various $y$ in $b$ represents a homology class in $L \times L$, which is related to $\tilde{\Phi}^1_{L,L}$ by Poincar{\'e} duality.
\end{remark}

We next turn to the multiplicative properties of \eqref{eq:corrected-phi}. By \eqref{eq:mu3} and \eqref{eq:phi2} we have
\begin{equation} \label{eq:derivation-property}
\begin{aligned}
& \tilde\phi^1_{L_0,L_2}(\mu^2(a_2,a_1)) - \mu^2(a_2,\tilde{\phi}^1_{L_0,L_1}(a_1)) - \mu^2(\tilde{\phi}^1_{L_1,L_2}(a_2),a_1) \\ & \qquad =
\mu^1\tilde{\phi}^2_{L_0,L_1,L_2}(a_2,a_1) + \tilde{\phi}^2_{L_0,L_1,L_2}(\mu^1(a_2),a_1)
+ (-1)^{|a_2|} \tilde{\phi}^2_{L_0,L_1,L_2}(a_2,\mu^1(a_1)), \\
& \tilde{\phi}^2_{L_0,L_1,L_2}(a_2,a_1) = \phi^2_{L_0,L_1,L_2}(b,a_2,a_1) \\ & \qquad \qquad
- \mu^3(c_{L_2},a_2,a_1) + \mu^3(a_2,c_{L_1},a_1) - \mu^3(a_2,a_1,c_{L_0}). \checked
\end{aligned}
\end{equation}
In particular, $\tilde\Phi^1_{L,L}$ is a derivation of the algebra $\HF^*(L,L)$, and therefore acts trivially on $\HF^0(L,L)$ (this could also be derived from Remark \ref{th:evaluation}, without using the multiplicative structure: $\tilde{\Phi}^1_{L,L}([e_L])$ is the represented by the evaluation map $(r,u) \longmapsto u(0,0)$, but that factors through the projection $S^1 \times \bar\scrM(y) \rightarrow \bar\scrM(y)$, hence represents the zero homology class.)

Suppose that we have an exact Lagrangian isotopy between $L_0$ and $L_1$, and choose gradings and {\it Spin} structures continuously over this isotopy. Isotopy invariance of Floer cohomology then provides a canonical cohomology class $[f] \in \HF^0(L_0,L_1)$, which is such that the multiplications
\begin{equation}
\CF^*(L_0,L_0) \xrightarrow{\mu^2(f,\cdot)} \CF^*(L_0,L_1)
\xleftarrow{\mu^2(\cdot,f)} \CF^*(L_1,L_1).
\end{equation}
are quasi-isomorphisms. Given a choice of cocycle $c_{L_0}$, one can always choose $c_{L_1}$ so that $\tilde{\Phi}^1_{L_0,L_1}$ is zero on $[f]$. This establishes a bijection between equivalence classes of $c_{L_0}$ and $c_{L_1}$. If we choose the $b$-equivariant structures in this way, the diagrams
\begin{equation}
\xymatrix{
\HF^*(L_1,L_2) \ar[d]_-{\tilde{\Phi}^1_{L_1,L_2}} \ar[r]^-{[\mu^2(\cdot,f)]} & \HF^*(L_0,L_2) \ar[d]^-{\tilde{\Phi}^1_{L_0,L_2}} \\
\HF^*(L_1,L_2) \ar[r]^-{[\mu^2(\cdot,f)]} & \HF^*(L_0,L_2)
}\qquad
\xymatrix{
\HF^*(L_2,L_0) \ar[d]_-{\tilde{\Phi}^1_{L_2,L_0}} \ar[r]^-{[\mu^2(f,\cdot)]} & \HF^*(L_2,L_1) \ar[d]_-{\tilde{\Phi}^1_{L_2,L_1}} \\ \HF^*(L_2,L_0) \ar[r]^-{[\mu^2(f,\cdot)]} & \HF^*(L_2,L_1)
}
\end{equation}
are commutative for any $b$-equivariant $L_2$. In particular,
\begin{equation}
L_1 \bullet_q L_2 = L_0 \bullet_q L_2, \quad L_2 \bullet_q L_1 = L_2 \bullet_q L_0.
\end{equation}
A parallel argument applies to changes of gradings instead of isotopies. If $L_1 = L_0[r]$, there is an element $[f]$ of degree $-r$ with the same properties. In particular,
\begin{equation} \label{eq:gr-shift}
L_0[r] \bullet_q L_2 = (-1)^r L_0 \bullet_q L_2, \quad
L_2 \bullet_q L_0[r] = (-1)^r L_2 \bullet_q L_0.
\end{equation}
(In contrast, changing the {\it Spin} structure affects $\phi^0_L$ in a more complicated way, and may not preserve $b$-compatibility in general.)

\begin{remark}
As a final comment on the various choices involved, suppose that we change $b$ to a cohomologous cocycle $b' = b + dy$. Given a choice of $y$, any $b$-equivariant Lagrangian submanifold becomes canonically $b'$-equivariant, by setting $c_L' = c_L + \phi^0_L(y)$. The difference this makes to the maps $\tilde\phi^1$ is, by \eqref{eq:phi1},
\begin{equation}
\phi^1_{L_0,L_1}(dy,a) + \mu^2(a,\phi^0_{L_0}(y)) - \mu^2(\phi^0_{L_1}(y),a) = - \mu^1(\phi^1_{L_0,L_1}(y,a)) - \phi^1_{L_0,L_1}(y,\mu^1(a)).
\checked
\end{equation}
Hence $\tilde{\Phi}^1$ and the $q$-intersection numbers remain unchanged.
\end{remark}

All we have said so far applies to arbitrary cocycles $b \in \CF^1(\lambda H)$. From this point onwards, we will use the dilation property \eqref{eq:liouville-element}, and choose a $\beta \in \CF^{-1}(\lambda H)$ such that
\begin{equation} \label{eq:corrected-liouville}
\delta b = e + d\beta.
\end{equation}
Then, \eqref{eq:hvee} and \eqref{eq:kvee} show that for any cocycle $a \in \CF^n(L,L)$,
\begin{equation}
\begin{aligned}
& e_L^\vee(\tilde{\phi}^1_{L,L}(a)) = e_L^\vee(\phi^1_{L,L}(b,a)) - e_L^\vee(\mu^2(c_L,a) - \mu^2(a,c_L)) \\
& = e_L^\vee(\phi^1_{L,L}(b,a)) - h^\vee_{L,L}(\phi_L^0(b),a) \\
& = \phi_L^{0,\vee}(\delta b,a) = e_L^\vee(a). \checked
\end{aligned}
\end{equation}
In the last step, we have used the fact that $\phi_L^{0,\vee}$ is a chain map, and that $\phi_L^{0,\vee}(e,\cdot)$ is chain homotopic to $e_L^\vee(\cdot)$. Hence, $\tilde\Phi^1_{L,L}$ acts as the identity on $\HF^n(L,L) \iso H^n(L;\K)$. As a consequence of this and the derivation property, one sees the following:

\begin{corollary} \label{th:poincare}
Under \eqref{eq:poincare}, the map $\tilde\Phi^1_{L_1,L_0}$ goes over to the dual (in the ordinary sense of linear algebra) of $\mathrm{id} - \tilde\Phi^1_{L_0,L_1}$. \qed
\end{corollary}

For the Euler characteristics this means that
\begin{equation}
L_1 \bullet_{q^{-1}} L_0 = (-1)^n q^{-1} (L_0 \bullet_q L_1).
\end{equation}

\begin{example}
Suppose that the algebra $H^*(L;\C)$ has a single generator of degree $n/k$. We know that $\tilde{\Phi}^1_{L,L}$ acts with weight $1$ on $H^n(L;\C)$. Since it is a derivation, it necessarily acts with weight $i/k$ on $H^{ni/k}(L;\C)$. In particular
\begin{equation}
L \bullet_q L = 1 + (-1)^{n/k} q^{1/k} + (-1)^{2 n/k} q^{2/k} + \cdots + (-1)^n q.
\end{equation}
\end{example}

\begin{remark} \label{th:theory}
The constructions in this section can be approached in a more conceptual way, provided that one is willing to take on board somewhat more homological algebra. Namely, let $\F(M)$ be the Fukaya $A_\infty$-category of $M$. There is a natural map
\begin{equation} \label{eq:full-open-closed}
\phi: \CF^*(\lambda H) \longrightarrow \CC^*(\F(M),\F(M))
\end{equation}
from the Hamiltonian Floer cochain complex to the Hochschild complex of that category \cite{seidel02,abouzaid10}, whose components extend the maps $\phi^k$ we defined for $k = 0,1,2$. From this and the given element $b$ one can form, in a purely algebraic way, a related category $\tilde{\F}(M,b)$ whose objects are $b$-equivariant Lagrangian submanifolds. The functor forgetting the equivariant structure maps $\tilde{\F}(M,b)$ fully and faithfully to the subcategory of $\F(M)$ consisting of those $L$ such that $\phi^0_L(b)$ is nullhomologous. $\tilde{\F}(M,b)$ carries its own Hochschild cocycle $\tilde{\phi}$, which is cohomologous to the pullback of $\phi(b)$ under the forgetful map, and such that the leading term $\tilde{\phi}^0$ is strictly zero. The next order terms $\tilde{\phi}^1, \tilde{\phi}^2, \dots$ then correspond to the objects of the same name we have defined. An alternative way to encode this cocycle would be to think of it as giving a first order deformation of $\tilde{\F}(M,b)$, parametrized by an infinitesimal parameter of degree $1$.
\end{remark}

\section{Dehn twists}\label{sec:dehn}

This section concerns one of the simplest classes of symplectic automorphisms, namely Dehn twists (also called Picard-Lefschetz monodromies) \cite{arnold95,seidel01}. We assume that $M$ has dimension $2n \geq 4$. This excludes the lowest-dimensional case of ordinary Dehn twists, which anyway would not be relevant in our context (see Example \ref{th:surfaces} below), and allows us to slightly simplify the discussion. Take a Lagrangian submanifold $V$ which is a sphere, and in fact comes with a diffeomorphism $S^n \rightarrow V$, determined up to isotopy and composition with $O(n+1)$. The associated Dehn twist $\tau_V$ has a natural grading, hence acts on objects of the Fukaya category.

Given $V$ and any two other Lagrangian submanifolds $(L_0,L_1)$, one considers the complex
\begin{equation}
\begin{aligned}
& T^*(L_0,L_1) = \CF^*(L_0,L_1) \oplus \Hom_\K^*(\CF(V,L_0),\CF(V,L_1))[-1], \\
& \mu^1_T(a,\alpha) = \big(\mu^1(a), v \mapsto -\mu^1(\alpha(v)) + (-1)^{|\alpha|}\alpha(\mu^1(v)) + \mu^2(a,v)\big). \checked
\end{aligned}
\end{equation}
Here, $|\alpha|$ is the natural degree of $\alpha$, before the shift has been applied. These complexes are covariantly functorial in $L_1$, by which we mean that they come with composition-type chain maps
\begin{equation} \label{eq:mu2-t}
\begin{aligned}
& \mu^2_T: \CF^*(L_1,L_2) \otimes T^*(L_0,L_1) \longrightarrow T^*(L_0,L_2), \\
& \mu^2_T(a_2,(a_1,\alpha_1)) = \big(\mu^2(a_2,a_1), v \mapsto (-1)^{|a_2|}\mu^2(a_2,\alpha_1(v))
- \mu^3(a_2,a_1,v)\big). \checked
\end{aligned}
\end{equation}
They are also contravariantly functorial in $L_0$, through maps
\begin{equation} \label{eq:mu2-tt}
\begin{aligned}
& \mu^2_T: T^*(L_1,L_2) \otimes \CF^*(L_0,L_1) \longrightarrow T^*(L_0,L_2), \\
& \mu^2_T((a_2,\alpha_2),a_1) = \big(\mu^2(a_2,a_1), v \mapsto \alpha_2(\mu^2(a_1,v)) - \mu^3(a_2,a_1,v)\big). \checked
\end{aligned}
\end{equation}
These maps have suitable homotopy associativity properties. We omit the general discussion, but we do want to mention one specific consequence. Namely, in \eqref{eq:mu2-t} suppose that $L_1 = L_2$. It is easy to see, by a suitable long exact sequence argument, that $\mu^2_T(e_{L_1},\cdot)$ induces an invertible map on the cohomology of $T^*(L_0,L_1)$. Homotopy associativity shows that this automorphism is idempotent, hence necessarily equal to the identity, thus showing that $\mu^2_T$ is homotopy unital, which is not immediately obvious since $e_{L_1}$ is not a strict unit. The same holds for \eqref{eq:mu2-tt}, of course.

From now on, we return to the usual situation where a dilation $b$ has been fixed. Because of our requirement that $n \geq 2$, $V$ itself can always be made $b$-equivariant. To simplify the following computations, let's assume that $\CF^*(V,V)$ is set up using a Hamiltonian $X_{\mu^1}$ which gives rise to the minimal number of chords. Then $\CF^1(V,V) = 0$, hence $\phi_V^0(b)$ will vanish strictly, and we can make $V$ equivariant by choosing $c_V = 0$. For any $L_0,L_1$ there is a map
\begin{equation}
\begin{aligned}
& \phi^1_T: T^*(L_0,L_1) \longrightarrow T^*(L_0,L_1), \\
& \phi^1_T(a,\alpha) = \big(\phi^1_{L_0,L_1}(b,a), v \mapsto \phi^1_{V,L_1}(b,\alpha(v))
- \alpha(\phi^1_{V,L_0}(b,v)) - \phi^2_{V,L_0,L_1}(b,a,v) \big).
\end{aligned}
\end{equation}
This satisfies the analogue of \eqref{eq:phi1}, meaning that
\begin{equation} \label{eq:phi1-t}
\mu^1_T(\phi^1_T(a,\alpha)) - \phi^1_T(\mu^1_T(a,\alpha)) = \mu^2_T(\phi^0_{L_1}(b),(a,\alpha)) - (-1)^{|a|} \mu^2_T((a,\alpha),\phi^0_{L_0}(b)). \checked
\end{equation}
Now suppose that $L_0,L_1$ are $b$-equivariant. Then one can add correction terms to make \eqref{eq:phi1-t} into a chain map
\begin{equation}
\begin{aligned}
& \tilde{\phi}^1_T: T^*(L_0,L_1) \longrightarrow T^*(L_0,L_1), \\
& \tilde{\phi}^1_T(a,\alpha)
= \phi^1_T(a,\alpha) - \mu^2_T(c_{L_1},(a,\alpha)) + \mu^2_T((a,\alpha),c_{L_0}) \\ &
\qquad = \big(\tilde\phi^1_{L_0,L_1}(a), v \mapsto \tilde\phi^1_{V,L_1}(\alpha(v)) - \alpha(\tilde\phi^1_{V,L_0}(v)) - \tilde\phi^2_{V,L_0,L_1}(a,v)\big), \checked
\end{aligned}
\end{equation}
where the last term is as in \eqref{eq:derivation-property}. This has a suitable derivation property: if $a_2 \in \CF^*(L_1,L_2)$ and $(a_1,\alpha_1) \in T^*(L_0,L_1)$ are cocycles, then
\begin{equation} \label{eq:tildephit-1}
\tilde\phi^1_T(\mu^2_T(a_2,(a_1,\alpha_1))) - \mu^2_T(a_2,\tilde{\phi}^1_T(a_1,\alpha_1)) - \mu^2(\tilde{\phi}^1_{L_1,L_2}(a_2),(a_1,\alpha_1)) = \mu^1_T(\text{something}).
\end{equation}
Similarly, on the other side we have
\begin{equation} \label{eq:tildephit-2}
\tilde\phi^1_T(\mu^2_T((a_2,\alpha_2),a_1)) - \mu^2_T((a_2,\alpha_2),\tilde{\phi}^1_{L_0,L_1}(a_1)) - \mu^2(\tilde{\phi}^1_T(a_2,\alpha_2),a_1) = \mu^1_T(\text{something}).
\end{equation}
The proof that this is the case is by a somewhat long computation, which we omit. It involves the next term $\mu^4$ in the Fukaya $A_\infty$-structure as well as the next order analogue $\phi^3$ of $\phi^2$.

So far, Dehn twists have not really played any role. All we need to know about such maps is contained in the following Lemma:

\begin{prop} \label{th:twist}
For every $L$ there is a cocycle $(x_L,\xi_L) \in T^0(L,\tau_V(L))$, canonical up to coboundaries, such that the following holds. Right composition yields quasi-isomorphisms
\begin{equation} \label{eq:quiso}
\mu^2_T(\cdot,(x_{L_0},\xi_{L_0})): \CF^*(\tau_V(L_0),L_1) \longrightarrow T^*(L_0,L_1).
\end{equation}
Similarly, left composition yields quasi-isomorphisms
\begin{equation}
\mu^2_T((x_{L_1},\xi_{L_1}),\cdot): \CF^*(L_0,L_1) \longrightarrow T^*(L_0,\tau_V(L_1)).
\end{equation}
Finally, the two maps fit into a homotopy commutative diagram
\begin{equation} \label{eq:lr-x}
\xymatrix{
\CF^*(L_0,\tau_V^{-1}(L_1)) \ar[dr]_-{\mu^2_T((x_{\tau_V^{-1}(L_1)},\xi_{\tau_V^{-1}(L_1)}),\cdot) \qquad\quad} \ar[rr]^-{(\tau_V)_*} && \CF^*(\tau_V(L_0),L_1) \ar[dl]^-{\qquad\mu^2_T(\cdot,(x_{L_0},\xi_{L_0}))} \\
& T^*(L_0,L_1) &
}
\end{equation}
where $(\tau_V)_*$ is the natural action of the symplectomorphism $\tau_V$ on Floer cochain complexes.
\end{prop}

\proof
The quasi-isomorphism properties are a consequence of \cite[Section 17]{seidel04}, which is based on  the geometric results from \cite{seidel01}. It remains to prove the homotopy commutativity of \eqref{eq:lr-x}, and to do that we have to recall some of the geometry involved. For any $L$, $x_L \in \CF^0(L,\tau_V(L))$ is obtained by counting pseudo-holomorphic sections of a Lefschetz fibration over the disc with one singular point, and which has $V$ as its vanishing cycle. We draw this schematically as in Figure \ref{fig:x} (the boundary marked point is a puncture, whereas the interior marked point is just the critical value of our Lefschetz fibration; the general style of these pictures is taken from \cite[Section 17]{seidel04}, to which we refer for a more detailed explanation). In parallel, Figure \ref{fig:x-homotopy} describes a one-parameter family of Lefschetz fibrations, whose parameter space can be thought of as a copy of the moduli space of discs with one interior and two boundary marked points (to simplify, one considers one interior and one boundary marked point as fixed, and lets the other one move). The associated parametrized section-counting invariant is a map $k: \CF^*(L_0,\tau_V^{-1}(L_1)) \rightarrow \CF^*(L_0,L_1)[-1]$ satisfying
\begin{equation}
\mu^1(k(a)) + k(\mu^1(a)) = \mu^2(x_{\tau_V^{-1}(L_1)},a) - \mu^2((\tau_V)_*(a),x_{L_0}).
\end{equation}%
\begin{figure}
\begin{centering}
\begin{picture}(0,0)%
\includegraphics{x.pstex}%
\end{picture}%
\setlength{\unitlength}{3355sp}%
\begingroup\makeatletter\ifx\SetFigFont\undefined%
\gdef\SetFigFont#1#2#3#4#5{%
  \reset@font\fontsize{#1}{#2pt}%
  \fontfamily{#3}\fontseries{#4}\fontshape{#5}%
  \selectfont}%
\fi\endgroup%
\begin{picture}(1560,2110)(703,-950)
\put(1426,-886){\makebox(0,0)[lb]{\smash{{\SetFigFont{10}{12.0}{\rmdefault}{\mddefault}{\updefault}{\color[rgb]{0,0,0}$L$}%
}}}}
\put(1351,989){\makebox(0,0)[lb]{\smash{{\SetFigFont{10}{12.0}{\rmdefault}{\mddefault}{\updefault}{\color[rgb]{0,0,0}$\tau_V(L)$}%
}}}}
\put(1601,389){\makebox(0,0)[lb]{\smash{{\SetFigFont{10}{12.0}{\rmdefault}{\mddefault}{\updefault}{\color[rgb]{0,0,0}$\tau_V$}%
}}}}
\end{picture}%
\caption{\label{fig:x}}
\end{centering}
\end{figure}%
\begin{figure}
\begin{centering}
\begin{picture}(0,0)%
\includegraphics{k.pstex}%
\end{picture}%
\setlength{\unitlength}{3355sp}%
\begingroup\makeatletter\ifx\SetFigFont\undefined%
\gdef\SetFigFont#1#2#3#4#5{%
  \reset@font\fontsize{#1}{#2pt}%
  \fontfamily{#3}\fontseries{#4}\fontshape{#5}%
  \selectfont}%
\fi\endgroup%
\begin{picture}(7343,2185)(-2264,-959)
\put(-1454,-536){\makebox(0,0)[lb]{\smash{{\SetFigFont{10}{12.0}{\rmdefault}{\mddefault}{\updefault}{\color[rgb]{0,0,0}$\tau_V^{-1}(L_1)$}%
}}}}
\put(-2249,-361){\makebox(0,0)[lb]{\smash{{\SetFigFont{10}{12.0}{\rmdefault}{\mddefault}{\updefault}{\color[rgb]{0,0,0}$L_0$}%
}}}}
\put(-1424,614){\makebox(0,0)[lb]{\smash{{\SetFigFont{10}{12.0}{\rmdefault}{\mddefault}{\updefault}{\color[rgb]{0,0,0}$L_1$}%
}}}}
\put(976,314){\makebox(0,0)[lb]{\smash{{\SetFigFont{10}{12.0}{\rmdefault}{\mddefault}{\updefault}{\color[rgb]{0,0,0}$\tau_V$}%
}}}}
\put(1351,-886){\makebox(0,0)[lb]{\smash{{\SetFigFont{10}{12.0}{\rmdefault}{\mddefault}{\updefault}{\color[rgb]{0,0,0}$L_0$}%
}}}}
\put(2026,614){\makebox(0,0)[lb]{\smash{{\SetFigFont{10}{12.0}{\rmdefault}{\mddefault}{\updefault}{\color[rgb]{0,0,0}$\tau_V^{-1}(L_1)$}%
}}}}
\put(626,614){\makebox(0,0)[lb]{\smash{{\SetFigFont{10}{12.0}{\rmdefault}{\mddefault}{\updefault}{\color[rgb]{0,0,0}$L_1$}%
}}}}
\put(2926,464){\makebox(0,0)[lb]{\smash{{\SetFigFont{10}{12.0}{\rmdefault}{\mddefault}{\updefault}{\color[rgb]{0,0,0}$L_1$}%
}}}}
\put(3781,619){\makebox(0,0)[lb]{\smash{{\SetFigFont{10}{12.0}{\rmdefault}{\mddefault}{\updefault}{\color[rgb]{0,0,0}$\tau_V(L_0)$}%
}}}}
\put(3901,-436){\makebox(0,0)[lb]{\smash{{\SetFigFont{10}{12.0}{\rmdefault}{\mddefault}{\updefault}{\color[rgb]{0,0,0}$L_0$}%
}}}}
\end{picture}%
\caption{\label{fig:x-homotopy}}
\end{centering}
\end{figure}%
Now, $\xi_L: \CF^*(V,L) \rightarrow \CF^*(V,\tau_V(L))[-1]$ itself is defined by counting sections in a similar one-parameter family, but where one of the two limits of the family contributes zero, by a vanishing result from \cite{seidel01}. This family and its limits are shown schematically in Figure \ref{fig:xi}. Finally, as indicated in Figure \ref{fig:kappa} one can define a two-parameter family which gives rise to a map $\kappa: \CF^*(L_0,\tau_V^{-1}(L_1)) \otimes \CF^*(V,L_0) \rightarrow CF^*(V,L_1)[-2]$,  satisfying
\begin{equation} \label{eq:k-kappa}
\begin{aligned}
& -\mu^1(\kappa(a_2,a_1)) + \kappa(\mu^1(a_2),a_1) + (-1)^{|a_2|} \kappa(a_2,\mu^1(a_1)) =
- \mu^2(k(a_2),a_1) \\ & + \mu^3((\tau_V)_*(a_2),x_{L_0},a_1) - (-1)^{|a_2|} \mu^2((\tau_V)_*(a_2),\xi_{L_0}(a_1)) \\ & + \xi_{\tau_V^{-1}(L)}(\mu^2(a_2,a_1))
- \mu^3(x_{\tau_V^{-1}(L_1)},a_2,a_1).
\end{aligned}
\end{equation}
The pair $(k,\kappa)$ yields precisely the homotopy needed in \eqref{eq:lr-x}. \qed
\begin{figure}
\begin{centering}
\begin{picture}(0,0)%
\includegraphics{xi.pstex}%
\end{picture}%
\setlength{\unitlength}{3355sp}%
\begingroup\makeatletter\ifx\SetFigFont\undefined%
\gdef\SetFigFont#1#2#3#4#5{%
  \reset@font\fontsize{#1}{#2pt}%
  \fontfamily{#3}\fontseries{#4}\fontshape{#5}%
  \selectfont}%
\fi\endgroup%
\begin{picture}(7343,2176)(-2264,-950)
\put(976,314){\makebox(0,0)[lb]{\smash{{\SetFigFont{10}{12.0}{\rmdefault}{\mddefault}{\updefault}{\color[rgb]{0,0,0}$\tau_V$}%
}}}}
\put(1351,-886){\makebox(0,0)[lb]{\smash{{\SetFigFont{10}{12.0}{\rmdefault}{\mddefault}{\updefault}{\color[rgb]{0,0,0}$V$}%
}}}}
\put(2026,614){\makebox(0,0)[lb]{\smash{{\SetFigFont{10}{12.0}{\rmdefault}{\mddefault}{\updefault}{\color[rgb]{0,0,0}$L$}%
}}}}
\put(376,614){\makebox(0,0)[lb]{\smash{{\SetFigFont{10}{12.0}{\rmdefault}{\mddefault}{\updefault}{\color[rgb]{0,0,0}$\tau_V(L)$}%
}}}}
\put(-2249,-361){\makebox(0,0)[lb]{\smash{{\SetFigFont{10}{12.0}{\rmdefault}{\mddefault}{\updefault}{\color[rgb]{0,0,0}$V$}%
}}}}
\put(2636,464){\makebox(0,0)[lb]{\smash{{\SetFigFont{10}{12.0}{\rmdefault}{\mddefault}{\updefault}{\color[rgb]{0,0,0}$\tau_V(L)$}%
}}}}
\put(-1274,-436){\makebox(0,0)[lb]{\smash{{\SetFigFont{10}{12.0}{\rmdefault}{\mddefault}{\updefault}{\color[rgb]{0,0,0}$L$}%
}}}}
\put(-1424,614){\makebox(0,0)[lb]{\smash{{\SetFigFont{10}{12.0}{\rmdefault}{\mddefault}{\updefault}{\color[rgb]{0,0,0}$\tau_V(L)$}%
}}}}
\put(3901,-436){\makebox(0,0)[lb]{\smash{{\SetFigFont{10}{12.0}{\rmdefault}{\mddefault}{\updefault}{\color[rgb]{0,0,0}$V$}%
}}}}
\put(3901,539){\makebox(0,0)[lb]{\smash{{\SetFigFont{10}{12.0}{\rmdefault}{\mddefault}{\updefault}{\color[rgb]{0,0,0}$V$}%
}}}}
\end{picture}%
\caption{\label{fig:xi}}
\end{centering}
\end{figure}%
\begin{figure}
\begin{centering}
\begin{picture}(0,0)%
\includegraphics{kappa.pstex}%
\end{picture}%
\setlength{\unitlength}{3355sp}%
\begingroup\makeatletter\ifx\SetFigFont\undefined%
\gdef\SetFigFont#1#2#3#4#5{%
  \reset@font\fontsize{#1}{#2pt}%
  \fontfamily{#3}\fontseries{#4}\fontshape{#5}%
  \selectfont}%
\fi\endgroup%
\begin{picture}(6096,3746)(-1547,-4435)
\put(1051,-2311){\makebox(0,0)[lb]{\smash{{\SetFigFont{10}{12.0}{\rmdefault}{\mddefault}{\updefault}{\color[rgb]{0,0,0}$\tau_V$}%
}}}}
\put(751,-3361){\makebox(0,0)[lb]{\smash{{\SetFigFont{10}{12.0}{\rmdefault}{\mddefault}{\updefault}{\color[rgb]{0,0,0}$V$}%
}}}}
\put(2251,-2311){\makebox(0,0)[lb]{\smash{{\SetFigFont{10}{12.0}{\rmdefault}{\mddefault}{\updefault}{\color[rgb]{0,0,0}$L_0$}%
}}}}
\put(526,-2236){\makebox(0,0)[lb]{\smash{{\SetFigFont{10}{12.0}{\rmdefault}{\mddefault}{\updefault}{\color[rgb]{0,0,0}$L_1$}%
}}}}
\end{picture}%
\caption{\label{fig:kappa}}
\end{centering}
\end{figure}

\begin{corollary}
$L$ admits an equivariant structure iff $\tau_V(L)$ does.
\end{corollary}

\proof From \eqref{eq:phi1-t} and the fact that $\mu^1_T(x_L,\xi_L) = 0$, we see that
\begin{equation}
\mu^2_T(\phi^0_{\tau_V(L)}(b),(x_L,\xi_L)) - \mu^2_T((x_L,\xi_L),\phi^0_L(b)) = \mu^1_T(\phi^1_T(x_L,\xi_L))
\end{equation}
is zero in cohomology. Multiplication with $(x_L,\xi_L)$ on either side is a quasi-isomorphism, and therefore, the vanishing of the class $[\phi^0_L(b)]$ is equivalent to that of $[\phi^0_{\tau_V(L)}(b)]$. In fact, \eqref{eq:lr-x} shows that more generally, these two obstruction classes are related by $(\tau_V)_*$. \qed

\begin{corollary} \label{th:induced-equivariant}
Given a $b$-equivariant structure on $L$, there is a unique such structure on $\tau_V(L)$ with the property that $\tilde\phi^1_T(x_L,\xi_L)$ is nullhomologous. Similarly, given a $b$-equivariant structure on $\tau_V(L)$, there is a unique such structure on $L$ with the corresponding property.
\end{corollary}

\proof $H^0(T(L,\tau_V(L))) \iso \HF^0(L,L)$ is one-dimensional. Changing the grading of $\tau_V(L)$ by a multiple of the identity element $e_{\tau_V(L)}$ changes the endomorphism $\tilde\phi_T^1$ of $T(L,\tau_V(L))$ by the corresponding multiple of $\mu^2_T(e_{\tau_V(L)},\cdot)$, which is homotopic to the identity by a previous observation. The statement follows directly from this, and the other part is parallel. \qed

\begin{corollary}
Suppose that the equivariant structure on $\tau_V(L_0)$ is induced from that on $L_0$ as in Corollary \ref{th:induced-equivariant}. Then this diagram is homotopy commutative:
\begin{equation} \label{eq:quiso-1}
\xymatrix{ \CF^*(\tau_V(L_0),L_1) \ar[rrrr]^-{\mu^2_T(\cdot,(x_{L_0},\xi_{L_0}))} \ar[d]_-{\tilde\phi^1_{\tau_V(L_0),L_1}} &&&& T^*(L_0,L_1) \ar[d]^-{\tilde\phi^1_T} \\
\CF^*(\tau_V(L_0),L_1) \ar[rrrr]^-{\mu^2_T(\cdot,(x_{L_0},\xi_{L_0}))} &&&& T^*(L_0,L_1).
}
\end{equation}
Correspondingly, suppose that the equivariant structure on $\tau_V^{-1}(L_1)$ is induced from that on $L_1$ in the same way. Then, the same conclusion holds for the following diagram:
\begin{equation} \label{eq:quiso-2}
\xymatrix{ \CF^*(L_0,\tau_V^{-1}(L_1)) \ar[d]_-{\tilde\phi^1_{L_0,\tau_V^{-1}(L_1)}} \ar[rrrr]^-{\mu^2_T((x_{\tau_V^{-1}(L_1)},\xi_{\tau_V^{-1}(L_1)}),\cdot)} &&&&
T^*(L_0,L_1) \ar[d]^-{\tilde{\phi}^1_T} \\
\CF^*(L_0,\tau_V^{-1}(L_1)) \ar[rrrr]^-{\mu^2_T((x_{\tau_V^{-1}(L_1)},\xi_{\tau_V^{-1}(L_1)}),\cdot)} &&&&
T^*(L_0,L_1).
}
\end{equation}
\end{corollary}

\proof Since we're talking about chain complexes of vector spaces, it is in fact sufficient to show that the diagrams commute on the cohomology level. But that follows directly from \eqref{eq:tildephit-1}, \eqref{eq:tildephit-2} and the fact that $\tilde{\phi}^1_T(x,\xi)$ is zero in cohomology. \qed

Commutativity of \eqref{eq:quiso-1} says that $\tilde{\Phi}^1_{\tau_V(L_0),L_1}$ is conjugate to the cohomology level map induced by $\tilde{\phi}^1_T$, which we write as $\tilde{\Phi}^1_T$. By construction, that maps sits in a diagram with long exact rows,
\begin{equation} \label{eq:les}
\xymatrix{
\cdots
H^*(T(L_0,L_1)) \ar[d]_-{\tilde{\Phi}^1_T}
\ar[r] & \HF^*(L_0,L_1) \ar[d]_-{\tilde{\Phi}^1_{L_0,L_1}} \ar[r] &
\Hom^*(\HF(V,L_0),\HF(V,L_1)) \ar[d] \cdots
 \\
\cdots H^*(T(L_0,L_1))
\ar[r] & \HF^*(L_0,L_1) \ar[r] &
\Hom^*(\HF(V,L_0),\HF(V,L_1)) \cdots
}
\end{equation}
where the right hand vertical map is
\begin{equation} \label{eq:tilde-end}
[\alpha] \longmapsto \tilde{\Phi}^1_{V,L_1} \circ [\alpha] - [\alpha] \circ \tilde{\Phi}^1_{V,L_0}.
\end{equation}
Suppose first that $L_0 = L_1 = V$. The right hand horizontal map in \eqref{eq:les} is dual to multiplication $[\mu^2]$. In this case, that map is injective, so $H^*(T(V,V)) \iso \HF^*(\tau_V(V),V)$ can be identified with a quotient of $\Hom_\K^*(H(V),H(V))[-1]$. More precisely, denote the standard generators of $H(V)$ by $[e]$ (identity) and $[f]$ (fundamental class), and their duals by $[e]^\vee$, $[f]^\vee$. The linear map \eqref{eq:tilde-end} acts on $[e] \otimes [e]^\vee$, $[e] \otimes [f]^\vee$, $[f] \otimes [e]^\vee$, $[f] \otimes [f]^\vee$, with eigenvalues $0$, $-1$, $1$, and $0$, respectively. If we divide out by the image of the dual of $[\mu^2]$, the outcome is that $\HF^*(\tau_V(V),V)$ has one generator in the two degrees $1-n$ and $1$, and that $\tilde{\Phi}^1$ acts on these with eigenvalues $-1$ and $0$, respectively. Since clearly $\tau_V(V) = V$ as a set, comparison with \eqref{eq:eq-shift}, \eqref{eq:gr-shift} shows the following:

\begin{corollary}
If one thinks of $\tau_V$ as acting on $b$-equivariant Lagrangian submanifolds in the way indicated by Corollary \ref{th:induced-equivariant}, then
\begin{equation}
\tau_V(V) = V[1-n]\langle 1 \rangle.
\qed
\end{equation}
\end{corollary}

Let's return to the case of a general $(L_1,L_2)$, and consider the implication for $q$-intersection numbers. Such numbers are additive in the long exact sequence \eqref{eq:les}, which is obvious from their formulation in terms of Lefschetz traces \eqref{eq:log-formula}. We can simplify \eqref{eq:tilde-end} by using Poincar{\'e} duality as in Corollary \ref{th:poincare}. The outcome is that $\tau_V(L_1) \bullet_q L_0 - L_1 \bullet_q L_0$ is, up to a sign $(-1)^{n+1}$, the supertrace of the endomorphism of $\HF^*(V,L_1) \otimes \HF^*(L_0,V)$ given by
\begin{equation}
\exp(\log(q) \tilde{\Phi}^1_{V,L_1}) \otimes \exp(\log(q) (\tilde{\Phi}^1_{L_0,V}-\id)).
\end{equation}
The second factor is just $q^{-1} \exp(\log(q) \tilde{\Phi}^1_{L_0,V})$. Hence we get the following {\em $q$-Picard-Lefschetz formula}:

\begin{theorem} \label{th:pl-q}
Let $(L_0,L_1)$ be $b$-equivariant Lagrangian submanifolds. For the induced $b$-equivariant structure on $\tau_V(L_0)$ as in Corollary \ref{th:induced-equivariant}, we have
\begin{equation} \label{eq:pl-q}
\tau_{V}(L_0) \bullet_q L_1 = L_0 \bullet_q L_1 + (-1)^{n+1} q^{-1} (L_0 \bullet_q V) (V \bullet_q L_1). \qed
\end{equation}
\end{theorem}

In general, even compactly supported symplectic automorphisms $\psi: M \rightarrow M$ may not preserve the class $[b] \in \HF^1(\lambda H)$, hence their action on Lagrangian Floer cohomology may not be compatible with the endomorphisms $\tilde{\Phi}^1$ (of course, this action should relate the endomorphisms obtained from $[b]$ to those for $\psi_*[b]$, but that is somewhat less useful in applications). However, for the special case of Dehn twists, we can show that this problem does not arise.

\begin{corollary}\label{cy:fun}
Let $(L_0,L_1)$ be $b$-equivariant Lagrangian submanifolds. For the induced $b$-equivariant structures on $(\tau_V(L_0),\tau_V(L_1))$, one has a commutative diagram
\begin{equation}
\xymatrix{
\HF^*(L_0,L_1) \ar[d]_-{\tilde{\Phi}^1_{L_0,L_1}}
\ar[rr]^-{(\tau_V)_*} && \HF^*(\tau_V(L_0),\tau_V(L_1))
\ar[d]^-{\tilde{\Phi}^1_{\tau_V(L_0),\tau_V(L_1)}} \\
\HF^*(L_0,L_1) \ar[rr]^-{(\tau_V)_*} && \HF^*(\tau_V(L_0),\tau_V(L_1)).
}
\end{equation}
\end{corollary}

\proof This follows directly by combining the two commutative diagrams \eqref{eq:quiso-1}, \eqref{eq:quiso-2} with \eqref{eq:lr-x}. \qed

\begin{remark}
We continue our discussion of the abstract framework from Remark \ref{th:theory}. To each $V$ one can associate an $A_\infty$-bimodule $\TT$ over $\F(M)$, whose underlying vector spaces are our $T(L_0,L_1)$, and whose structure maps include our $\mu^1_T$, $\mu^2_T$. Take this bimodule and pull it back to $\tilde{\F}(M,b)$ via the forgetful map. The algebraic part of our argument says that this admits a canonical extension over the infinitesimal deformation of $\tilde{\F}(M,b)$ determined by $\tilde{\phi}$. The map $\tilde\phi^1_T$ is the deformed differential, which is the simplest term of that structure.

The relation between $\TT$ and the geometric Dehn twist $\tau_V$ is as follows. Pull back $\TT$ by $\tau_V: \F(M) \rightarrow \F(M)$, acting on the left. The result should then be quasi-isomorphic to the diagonal bimodule, with the quasi-isomorphism given by a cocycle in the Hochschild complex of $\F(M)$ with coefficients in $(\tau_V)^*\TT$. The data provided in Proposition \ref{th:twist} are the first two terms of such a cocycle. There is no description of the higher order terms in the existing literature, as far as the authors know, even though an equivalent statement should follow from ongoing work of Ma'u-Wehrheim-Woodward \cite{mww,ww}.
\end{remark}

\section{Liouville manifolds\label{sec:viterbo}}

We want to specialize Assumption \ref{th:convexity} further, in order to arrive at the situation considered in the Introduction. Suppose that $(\bar{M}, \bar\omega = d\bar\theta)$ is a Liouville domain. This means that $\bar{M}$ is a compact manifold with boundary, with an exact symplectic form, and such that the Liouville vector field dual to $\bar\theta$ points strictly outwards along the boundary. We also assume that this comes with a distinguished homotopy class of trivializations of its canonical bundle. There is a preferred way to attach an infinite cone to the boundary, forming a Liouville manifold (a particular kind of noncompact symplectic manifold)
\begin{equation}
M = \bar{M} \cup_{\partial \bar{M}} (\R^+ \times \partial \bar{M}).
\end{equation}
$M$ carries an exact symplectic structure $\omega = d\theta$ which extends the given one on $\bar{M}$, and such that $\theta = e^r(\bar\theta|\partial \bar{M})$ on the cone part. There is a standard class of Hamiltonian functions on $M$, namely those which satisfy $H(r,y) = e^r$ on the cone. Likewise, there is a standard class of compatible almost complex structures $I$, namely those whose restriction to the cone is translation-invariant in $r$-direction and satisfies $d(e^r) \circ I = -\theta$. If we pick any such $I$ and $H$, then Assumption \ref{th:convexity} will be satisfied for the exhaustion by large level sets of $H$. This is a standard maximum principle argument.

Taking $H$ within this class, we can consider its Floer cohomology $\HF^*(\lambda H)$. For sufficiently small $\lambda>0$, this is isomorphic to the ordinary Floer cohomology $H^*(M;\K)$ \cite{viterbo97a}. On the other hand, we can pass to the direct limit of the maps \eqref{eq:continuation} as $\lambda \rightarrow \infty$. The resulting graded vector space
\begin{equation}
\SH^*(M) \stackrel{\mathrm{def}}{=} \underrightarrow{\lim}_{\lambda} \, \HF^*(\lambda H)
\end{equation}
is a symplectic invariant of $M$, called {\em symplectic cohomology} (with orientations twisted by the class $[\alpha] \in H^2(M;\Z/2)$, to be precise). The original reference is \cite{viterbo97a}, see also \cite{cieliebak-floer-hofer95} for a different approach; more recent surveys are \cite{oancea04b,seidel07}. By construction, it comes with a map $H^*(M;\K) \rightarrow \SH^*(M)$. One can show that the image of the ordinary unit in cohomology under this map coincides with the limit of the canonical elements $[e] \in \HF^0(\lambda H)$. Let's denote this image class by $E \in \SH^0(M)$. Similarly, the BV operators induce a map $\Delta: \SH^*(M) \rightarrow \SH^{*-1}(M)$, which vanishes on the image of $H^*(M;\K) \rightarrow \SH^*(M)$. With this in mind, we return to the main thread of our discussion. By construction, the existence of a dilation $b \in \CF^1(\lambda H)$ for some $\lambda > 0$ implies that we have a class $B \in \SH^1(M)$ such that
\begin{equation} \label{eq:Delta}
\Delta B = E.
\end{equation}
Conversely, if there is such a class, then dilations exist at least for $\lambda \gg 0$. By a slight abuse of terminology, we will refer to a solution $B$ of \eqref{eq:Delta} also as a dilation. Such elements will be the main focus of the following discussion.

\begin{example} \label{th:surfaces}
In the lowest dimension $2n = 2$, the situation is as follows. If $M = \R^2$, then $\SH^*(M) = 0$, so $B = 0$ is a trivial solution of \eqref{eq:Delta}. In all other cases, $\SH^*(M)$ splits into $H^*(M;\K)$ and another summand measuring the contribution of non-contractible periodic orbits, see for instance \cite[Section 8]{bourgeois-oancea09}. Moreover, this splitting is compatible with $\Delta$, which implies that we cannot have any dilations.
\end{example}

Symplectic cohomology has been most thoroughly studied for cotangent bundles $M = T^*L$ of closed oriented manifolds $L$. Take $[\alpha] \in H^2(M;\Z/2)$ to be the pullback of $w_2(L)$. Then, one has a canonical isomorphism with free loop space homology, 
\begin{equation} \label{eq:loop}
\SH^*(M) \iso H_{n-*}(\scrL L;\K).
\end{equation}
The original proof is in \cite{viterbo97b}, and two others have appeared since then \cite{salamon-weber03, abbondandolo-schwarz06} (the comparison issue for coherent orientations, which requires this particular choice of $[\alpha]$, has been clarified in \cite{kragh07,seidel10,abouzaid10b}). In this context, the map $H^*(M;\K) \rightarrow \SH^*(M)$ is just the inclusion of constant loops. Moreover, the BV operator corresponds to the map
\begin{equation} \label{eq:rotation}
H_*(\scrL L;\K) \hookrightarrow H_{*+1}(S^1 \times \scrL L;\K) \longrightarrow H_{*+1}(\scrL L;\K),
\end{equation}
where the second arrow is rotation of loops (the last-mentioned fact follows from a suitable Morse-Bott setup and the arguments in \cite{abbondandolo-schwarz06}, see the discussion in \cite[Introduction]{bourgeois-oancea09b}).

\begin{example} \label{th:cotangent}
If $L$ is a $K(\pi,1)$ space, the situation is quite similar to the case of noncontractible surfaces, in that we have a topological splitting of $\SH^*(T^*L)$ which excludes the existence of dilations. This conclusion holds for arbitrary $[\alpha]$, using a mild generalization of \eqref{eq:loop}. 
\end{example}

This observation has a wider implication, which severely restricts the class of Liouville manifolds which can potentially admit dilations.

\begin{corollary}
Suppose that $M$ is a Liouville manifold admitting a dilation. Then, it cannot contain a closed exact Lagrangian submanifold $L \subset M$ which is a $K(\pi,1)$ space.
\end{corollary}

\proof This is a direct consequence of Example \ref{th:cotangent} and Viterbo functoriality \cite{viterbo97b}, which for any exact $L \subset M$ provides a map $\SH^*(M) \rightarrow \SH^*(T^*L)$ compatible with the identity elements and dilations. In fact, the conclusion is somewhat similar to \cite[Corollary 5.4]{viterbo97a}. \qed

\begin{example}
Still within the class of cotangent bundles $M = T^*L$, take $L = S^n$ for $n \geq 2$. The Hurewicz isomorphism $\pi_n(L) \rightarrow H_n(L; \Z)$ factors through
\begin{equation}
\pi_n(L) = \pi_{n-1}(\Omega L) \longrightarrow H_{n-1}(\scrL L;\Z) \longrightarrow H_n(\scrL L;\Z) \longrightarrow H_n(L;\Z),
\end{equation}
where the middle arrow is \eqref{eq:rotation}, and the right one is evaluation of loops at endpoints. For $n>2$ the last-mentioned map is an isomorphism; whereas for $n = 2$ it has a $\Z/2$ kernel \cite[Theorem 2]{cohen-jones-yan02}. Hence, a solution of \eqref{eq:Delta} exists for arbitrary $\K$ if $n>2$, and for $\mathrm{char}(\K) \neq 2$ if $n = 2$ (a delicate computation from \cite{menichi06} shows that the assumption on the characteristic cannot be removed).

More generally, suppose that $L$ admits a nonzero degree map $S^n \times \tilde{L} \rightarrow L$, for some $n \geq 2$ and arbitrary $\tilde{L}$. Let's take our coefficient field $\K$ to be of characteristic zero. Then $T^*L$ (equipped with $[\alpha] = w_2(L)$) again admits a dilation, inherited from that on the sphere. This applies for instance to complex projective spaces, since those spaces are the symmetric products of the two-sphere (there are also slightly more direct ways of verifying the existence of dilations on $T^*\C P^{n/2}$, for instance using rational homotopy theory \cite{burghelea-vigue-poirrier85}). 
\end{example}

\begin{remark}
Symplectic cohomology has an equivariant version, defined as
\begin{equation}
\SH^*_{\eq}(M) = \underrightarrow{\lim}_{\lambda} \, \HF^*_{\eq}(\lambda H)
\end{equation}
where $\HF^*_{\eq}(\lambda H)$ is as in Remark \ref{th:equivariant}, and the maps are obtained by adding suitable higher order correction terms to \eqref{eq:continuation}. For small $\lambda>0$ we again have $\HF^*_{\eq}(\lambda H) \iso H^*(M;\K[u^{-1}])$, which gives rise to a map
\begin{equation} \label{eq:eq-accelleration}
H^*(M;\K[u^{-1}]) \longrightarrow \SH^*_{\eq}(M).
\end{equation}
Now suppose that $\CF^*(\lambda H)$ contains a dilation $b$, which means that $db = 0$, $\delta b = e + d\beta$ for some $\beta$. By definition \eqref{eq:equivariant-differential}, this implies that $d_{\mathrm{eq}}(-\beta + u^{-1}b) = e$. Hence, if \eqref{eq:Delta} has a solution, the image of $1 \in H^0(M;\K)$ under \eqref{eq:eq-accelleration} must vanish. However, there is no reason to expect the converse to be true as well.

Equivariant symplectic cohomology is particularly well suited for making the connection with symplectic field theory (SFT). Work in progress by Bourgeois-Oancea, following on from \cite{bourgeois-oancea09}, will show that for $\mathrm{char}(\K) = 0$, \eqref{eq:eq-accelleration} fits into a long exact sequence
\begin{equation} \label{eq:oancea}
\cdots \rightarrow \mathit{CH}_{2n-2-*}(\partial \bar{M}) \longrightarrow H^*(M;\K[u^{-1}]) \longrightarrow \SH^*_{\eq}(M) \rightarrow \cdots
\end{equation}
where $\mathit{CH}_*(\partial \bar{M})$ is the linearized contact homology of the boundary of the Liouville domain $\bar{M}$ from which $M$ is constructed. The first map in \eqref{eq:oancea} is defined in terms of relative SFT invariants. This is particularly illuminating in the algebro-geometric context, where $M = X \setminus D$ is the complement of a smooth ample divisor $D$ inside a projective algebraic variety $X$. One then expects that the SFT invariants can also be expressed in terms of relative Gromov-Witten invariants of $(X,D)$ (this has not been proved, as far as we know, but see for instance \cite{katz07} for a comparison of the two formalisms).
\end{remark}

\section{Lefschetz fibrations}\label{sec:lefschetz}
Remaining in the framework of Section \ref{sec:viterbo}, it is natural to want to go beyond the case of cotangent bundles. In order to do that, one needs some way of making at least partial computations in symplectic cohomology. Here, we use McLean's approach via Lefschetz fibrations \cite{mclean09} to provide some additional examples where dilations can be shown to exist.

Working with Lefschetz fibrations which have noncompact fibres requires some technical caution. We start by describing the desired structure of such a fibration at infinity, which in particular means away from the singularities. The {\em base} will be $B = \C$ with its standard exact symplectic structure $\omega_B = d\theta_B$, $\theta_B = \frac{i}{4}(z\, d\bar{z} - \bar{z}\, dz)$, and complex structure $I_B = i$. The {\em fibre} will be a $(2n-2)$-dimensional Liouville manifold, constructed as $F = \bar{F} \cup_{\partial\bar{F}} (\R^+ \times \partial \bar{F})$, and with a fixed almost complex structure $I_F$ which is of the same type as in Section \ref{sec:viterbo}. The {\em outer monodromy} $\mu$ is a symplectic automorphism of $F$, which is the identity outside $\bar{F}$, and such that $\mu^*\theta_F - \theta_F$ is the derivative of a function $R$ supported in $\bar{F}$. Let $T = (\R \times F)/(t,x) \sim (t-1,\mu(x))$ be its mapping torus, equipped with the one-form $\theta_T = \theta_F + d(tR)$. Define
\begin{equation}
\tilde{M} = (\R^+ \times T) \cup_{(\R^+ \times S^1 \times \R^+ \times \partial\bar{F})} (\C \times \R^+ \times \partial \bar{F}),
\end{equation}
where the two parts are joined together by identifying $(s,t,r,y)$ with $(\exp(2\pi(s+it)),r,y)$. This comes with a map $\tilde\pi: \tilde{M} \rightarrow \C$, which is $\tilde\pi(s,t,x) = \exp(2\pi(s+it))$ and $\tilde\pi(z,r,y) = z$ on the two parts. Outside the unit disc, this is a fibration with fibre $F$ and monodromy $\mu$. Inside the unit disc, the fibres consist only of the conical part of $F$, and the fibration is trivial there. Equip $\tilde{M}$ with a one-form $\theta_{\tilde{M}}$ which equals  $\theta_T + \tilde\pi^*\theta_B$ and $\theta_F + \tilde\pi^*\theta_B$ on the two pieces, and with its exterior derivative $\omega_{\tilde{M}} = d\theta_{\tilde{M}}$. Take a one-parameter family $I_s$ of compatible almost complex structures on $F$ which all agree with $I_F$ on the conical part, and satisfy $\mu^*I_{s-1} = I_s$. This can be viewed as a family of almost complex structures on the fibres of $T \rightarrow S^1$, and it combines with $I_B$ to yield an almost complex structure on $\R^+ \times T$. Extend that by the product structure $I_B \times I_F$ over the rest of $\tilde{M}$, and denote the outcome by $I_{\tilde{M}}$.

\begin{definition} \label{def:lefschetz}
A Lefschetz fibration $\pi: M \rightarrow B$, with fibre $F$ and outer monodromy $\mu$, is given by the following data. The total space is an exact symplectic manifold $(M,\omega_M = d\theta_M)$, equipped with a compatible almost complex structure $I_M$. The map $\pi$ is $I_M$-holomorphic and has only finitely many critical points, which are locally modelled on nondegenerate critical points of holomorphic functions. There is a relatively compact open subset of $M$ whose complement can be identified with the $\tilde{M}$ constructed above, in a way which is compatible with all the additional structure. Finally, we want to have a preferred trivialization of the canonical bundle $\scrK_M$ (this induces the same structure on the fibre $F$, at least up to homotopy).
\end{definition}

Given such a Lefschetz fibration, choose a Hamiltonian function $H_F$ on $F$ which vanishes outside the conical end, and which on that end is of the form $H_F(r,x) = \psi(e^r)$, where $\psi(t) = 0$ for $t \leq 1$, $\psi''(t) = 0$ for $t \geq 2$, and $\psi''(t) > 0$ for $t \in (1,2)$. We assume that there are no Reeb orbits of length $\lambda = \psi'(1)$ on $\partial\bar{F}$, so that $\HF^*(H_F)$ is well-defined. Choose also $b > 1$, and take the function $H_B(z) = \epsilon |z-b|^2/2$ on the base, for some small $\epsilon>0$. These two combine to yield a function on the total space, schematically denoted by
\begin{equation} \label{eq:hm-definition}
H_M = H_B + H_F.
\end{equation}
More precisely, we take $H_F$ on the conical part of each fibre, where the fibration is trivial, extend it by zero over the rest, and then add the pullback of $H_B$. Outside the preimage of the unit disc, the Hamiltonian vector field of $H_M$ projects to that of $H_B$. Using this, one can show that Assumption \ref{th:convexity} is satisfied, and that $\HF^*(H_B)$ is defined.

\begin{lemma} \label{th:lefschetz-sequence}
There is a long exact sequence
\begin{equation} \label{eq:lefschetz-sequence}
\cdots \rightarrow \K^{\mathit{Crit}(\pi)}[-n] \longrightarrow \HF^*(H_M) \longrightarrow \HF^*(H_F) \rightarrow \cdots
\end{equation}
where the left hand term has one copy of $\K$ for each critical point of $\pi$, all placed in degree $n = \mathrm{dim}_\C(M)$.
\end{lemma}

\proof Even though the definitions involve perturbations of the Hamiltonian and almost complex structures, it is useful to first look at the unperturbed situation.

The one-periodic orbits of the Hamiltonian vector field $X_M$ associated to $H_M$ are of three kinds:
(i) if we identify the fibre $\pi^{-1}(b)$ with $F$ in the obvious way, then $H_M|\pi^{-1}(b) = H_F$ vanishes in $\bar{F}$. All points of that subset are critical points of $H_M$, hence constant one-periodic orbits. The value of the action functional \eqref{eq:action} on those points is zero. (ii) still within the same fibre, each periodic orbit of the Reeb flow on $\partial\bar{F}$ with period $0 < l < \lambda$ gives rise to a one-periodic orbit of $H_M$. The value of the action functional on such points is $\psi(e^r) - e^r\psi'(e^r)$, where $r$ is the unique point such that $\psi'(e^r) = l$. By assumption on $\psi$, these values are negative (this is a familiar argument in symplectic homology theory, compare for instance \cite[Section 3c]{seidel07}). (iii) Any critical point of $\pi$ gives rise to a constant one-periodic orbit, with action value $\epsilon |\pi(x) - b|^2/2 > 0$. Moreover, these orbits are nondegenerate, and their Conley-Zehnder index equals their Morse index, which is $n$.

Look at a solution of the unperturbed Floer equation with limits of type (i) or (ii), meaning
\begin{equation} \label{eq:bare-floer}
\left\{
\begin{aligned}
& u: \R \times S^1 \longrightarrow M, \\
& \partial_s u + I_M(\partial_t u - X_M) = 0, \\
& \textstyle\lim_{s \rightarrow \pm\infty} u(s,t) = y_\pm(t) \in \pi^{-1}(b).
\end{aligned}
\right.
\end{equation}
We will show that then, the whole of $u$ is necessarily contained in the same fibre. By construction, at all points where $u(s,t)$ lies sufficiently close to that fibre, $v = \pi(u)$ satisfies the equation
\begin{equation} \label{eq:base-j}
\partial_s v + i(\partial_t v - X_B) = 0,
\end{equation}
where $X_B(z) = \epsilon i (z-b)$ is the infinitesimal rotation vector field centered at $b$. Take a two-form $\nu \in \Omega^2(B)$ which is supported in a small neighbourhood of $b$, is positive at that point and nonnegative everywhere, and moreover is invariant under $X_B$. One can then write $\nu(\cdot,X_B) = dN$ for some function $N$. Given a solution $u$ of \eqref{eq:bare-floer}, we can compute
\begin{equation} \label{eq:nu-argument}
\int_{\R \times S^1} \nu(\partial_s v, \partial_t v - X_B) ds \wedge dt = \int_{\R \times S^1} v^*\nu - dN(\partial_s v) ds \wedge dt
\end{equation}
in two ways. On one hand, by \eqref{eq:base-j} the integrand is equal to $\nu(\partial_s v, i\partial_s v)$, hence nonnegative everywhere. On the other hand, applying Stokes shows that the overall integral is zero. Comparing the two ideas shows that $\partial_s v$ must be zero at all points where $v(s,t)$ lies within the support of $\nu$. Hence, $v$ is necessarily constant equal to $b$, which in turn implies the desired statement that $u$ remains inside $\pi^{-1}(b)$.

When we perturb the Hamiltonian, that can be done so that near $\pi^{-1}(b)$ it equals $H_B$ plus some time-dependent function on the fibre. We then have one-periodic orbits contained inside $\pi^{-1}(b)$, which arise from those of types (i) and (ii) mentioned above. Moreover, each point of type (iii) will give rise to a nearby one-periodic orbit with the same properties. Finally, the argument about solutions $u$ goes through without any changes. The upshot is that the generators of $\CF^*(H_M)$ lying close to the critical points of $\pi$ form a subcomplex, whose differential is necessarily trivial for degree reasons; and moreover, the quotient complex can be identified with $\CF^*(H_F)$. \eqref{eq:lefschetz-sequence} is the induced long exact sequence. \qed

An argument of the same kind shows that the map $\HF^*(H_M) \rightarrow \HF^*(H_F)$ from \eqref{eq:lefschetz-sequence} is compatible with the basic additional structures, meaning the identity elements and the BV operators. The map $HF^0(H_M) \rightarrow HF^0(H_F)$ is always injective since by definition $n \geq 1.$ Next, note that in the construction $\lambda$ can be chosen as large as desired. Hence, if $\SH^*(F)$ admits a solution of \eqref{eq:Delta}, then so does $\HF^*(H_M)$ for suitably chosen $H_M$, provided that $n>2$, which ensures that $\HF^1(H_M) \rightarrow \HF^1(H_F)$ is onto. Finally, the total space $M$ is itself again a Liouville manifold, and it is not difficult to construct a continuation map $\HF^*(H_M) \rightarrow \SH^*(M)$, with the same compatibility properties as before. We omit the details, referring to \cite{oancea08,mclean09} for similar arguments, but the outcome is the following:

\begin{proposition} \label{th:lefschetz}
If $\SH^*(F)$ contains a dilation, and $n>2$, then so does $\SH^*(M)$. \qed
\end{proposition}

\begin{example} \label{th:an}
Consider the Milnor fibre of the $n$-dimensional $(A_m)$ singularity for some $m \geq 1$, $n \geq 3$:
\begin{equation} \label{eq:milnor}
M = \{x \in \C^4 \;:\; p(x) = x_1^{m+1} - 1 + x_2^2 + \cdots + x_{n+1}^2 = 0\}.
\end{equation}
Projection to the $x_1$-variable turns this into a Lefschetz fibration whose smooth fibres are regular affine conics, hence symplectically isomorphic to $T^*S^{n-1}$. Granted, this does not satisfy the triviality condition from Definition \ref{def:lefschetz}, but that can be repaired by suitably deforming the symplectic form and almost complex structure, compare \cite[Section 19]{seidel04}. Since $M$ is a hypersurface and $H^1(M) = H^2(M) = 0$, there is a unique trivialization of $\scrK_M$ up to homotopy, and a unique choice of $[\alpha]$. By combining Example \ref{th:cotangent} and Proposition \ref{th:lefschetz}, one finds that $M$ admits a solution of \eqref{eq:Delta} (it is plausible that the same should be true for $n = 2$, but there is no proof of this at present, except in the trivial case $m = 1$). 

From now on, we specialize to $n = 3$, and take $\mathrm{char}(\K) = 0$. The symplectic geometry of \eqref{eq:milnor} has been studied repeatedly, see for instance \cite{khovanov-seidel98}. Given an embedded path $\gamma \subset \C$ whose endpoints are $(m+1)$-st roots of unity, and which does not otherwise intersect any $(m+1)$-st roots of unity, there is an associated Lagrangian three-sphere $L_{\gamma} \subset M$. In particular, given a chain of paths $(\gamma_1,\dots,\gamma_m)$ which forms an $(A_m)$ type diagram, one gets a corresponding configuration of Lagrangian spheres $L_i = L_{\gamma_i}$. One can adjust the gradings of these spheres so that
\begin{equation} \label{eq:neighbours}
\HF^k(L_i,L_{i+1}) = \begin{cases} \K & k = 1, \\ 0 & \text{in all other degrees.} \end{cases}
\end{equation}
Fix a dilation $b$, and make the $L_i$ into $b$-equivariant Lagrangian submanifolds, in such a way that $\tilde\Phi^1$ acts as $1/3$ on \eqref{eq:neighbours}. Then, as a consequence of the derivation property and Poincar{\'e} duality, all endomorphisms $\tilde\Phi^1$ multiply the degree $k$ part of $\mathit{HF}^*(L_i,L_j)$ by $k/3$.

Since any $L_{\gamma}$ can be obtained from our basic collection $(L_1,\dots,L_m)$ by applying Dehn twists, their $q$-intersection numbers can be computed, at least up to powers of $q$, from the given one by applying Proposition \ref{th:pl-q}. The answer turns out to be equivalent to improved intersection numbers introduced by Givental in \cite{givental88} (in spite of this, there are obvious philosophical differences between the two approaches; Givental's definition of $q$-intersection numbers uses the entire Milnor fibration, rather than just a fibre, and is also highly unsymmetric).
\end{example}


For the next example, we need a slight variation of the previous discussion. Suppose that we have a Lefschetz fibration $\pi: M \rightarrow B$ whose base is a once-punctured plane (say $B = \C/i\Z$, with symplectic form $\frac{i}{2} dz \wedge d\bar{z}$). There are now outer monodromies around both $\mathrm{re}(z) = -\infty$ and $\mathrm{re}(z) = +\infty$, and we assume that the first of these is trivial.
Otherwise, the assumptions are analogous to those in Definition \ref{def:lefschetz}. On the base choose the function $H_B(z) = \epsilon(\mathrm{re}(z) - b)^2/2$, where the constants are $b \ll 0$ and a small $\epsilon>0$. Use this and a function $H_F$ on the fibre to define $H_M$ as in \eqref{eq:hm-definition}. Then, an argument similar to the proof of Lemma \ref{th:lefschetz-sequence} applies, which we will now sketch briefly. Critical points of type (i) form a copy of $\bar{F} \times S^1$, and have action value $0$. There is an $S^1$ of critical points of type (ii) for each one-periodic point of the flow of $H_F$, with negative action value. Finally, critical points of type (iii) again come from the singularities of the Lefschetz fibration, and have positive actions. In the analogue of \eqref{eq:base-j}, we have $X_B(z) = \epsilon i (\mathrm{re}(z) - b)$. Define $\nu$ to be supported in a small neighbourhood of $\{\mathrm{re}(w) = b\}$ and invariant under $X_B$. Then, the argument concerning \eqref{eq:nu-argument} applies as before, allowing us to determine the subcomplex consisting of generators of type (ii) and (iii). The outcome is a long exact sequence analogous to \eqref{eq:lefschetz-sequence}, of the form
\begin{equation} \label{eq:lefschetz-sequence-2}
\cdots \rightarrow \K^{\mathit{Crit}(\pi)}[-n] \longrightarrow \HF^*(H_M) \longrightarrow \HF^*(H_F) \otimes H^*(S^1;\K) \rightarrow \cdots
\end{equation}
(for another occurrence of the same geometric situation, see \cite[Section 3]{albers-mclean09}). By the same kind of argument, the BV operator on the last term is $[\delta] \otimes \mathrm{id}$. This leads to the same consequence as in Proposition \ref{th:lefschetz}.

\begin{example}
Consider what one might call the three-dimensional ``affine $(A_1)$ type Milnor fibre'', namely
\begin{equation}
M = \{x \in \C^* \times \C^3 \;:\; p(x) = x_1^2 - 1 + x_2^2 + x_3^2 + x_4^2 = 0\}
\end{equation}
with the trivialization of its canonical bundle given by $(dx_1/x_1 \wedge dx_2 \wedge dx_3 \wedge dx_4)/dp$. The geometry is similar to the $m = 1$ case of Example \ref{th:an}, except that Lagrangian spheres $L_\gamma \subset M$ are now associated to paths $\gamma$ in $\C^*$ with endpoints $\pm 1$. Take the two spheres $L_0,L_1$ associated to arcs in the lower and upper half-planes, respectively. For an appropriate choice of gradings,
\begin{equation}
\HF^k(L_0,L_1) = \begin{cases} \K & k = 1,2, \\ 0 & \text{in all other degrees.} \end{cases}
\end{equation}
We know from the generalization of Proposition \ref{th:lefschetz} mentioned above that $M$ admits a solution of \eqref{eq:Delta}. In fact, again taking $\K = \Q$ for simplicity, it admits an infinite family of such solutions, of the form $B + \mu Z$, where $Z$ is the image in $\SH^1(M)$ of the generator of $H^1(M)$. Passing from $B$ to $B + \mu Z$ changes the action of $\tilde\Phi^1$ on $\HF^1(L_0,L_1)$ by some $\lambda$, and that on $\HF^2(L_0,L_1)$ by $\lambda + \mu$. Hence, by a suitable choice of dilation $b$ and of the $b$-equivariant structure, one can achieve that $\tilde{\Phi}^1$ acts as $k/3$ on $\HF^k(L_0,L_1)$. The same will then be true for the other groups $\HF^*(L_i,L_j)$ as well.
\end{example}

\begin{remark}
In the examples above, we have paid particular attention to the case when the action of $\tilde{\Phi}^1$ is proportional to the grading. This is motivated by the following purely algebraic observation. Let $\A$ be an $A_\infty$-algebra over a field $\K$ of characteristic zero, and $A = H(\A)$ its cohomology algebra. There are canonical maps between Hochschild cohomology groups,
\begin{equation}
\begin{aligned}
& \Psi^0: \HH^*(\A,\A) \longrightarrow \HH^0(A,A), \\
& \Psi^1: \mathit{ker}(\Phi^0) \longrightarrow \HH^1(A,A)[-1].
\end{aligned}
\end{equation}
Here, $\HH^0(A,A) = Z(A)$ is the center of $A$, which is itself graded, and similarly for $\HH^1(A,A) = \mathit{Der}(A,A)/\mathit{Inn}(A,A)$. Suppose that there is a $B \in \HH^1(\A,\A)$ such that $\Psi^0(B) = 0$, and $\Psi^1(B)$ agrees with the Euler derivation (the one that multiplies the degree $i$ part of $A$ by $i$) up to inner derivations. One can prove that if this is the case, then $\A$ is formal (quasi-isomorphic to its cohomology). Via the abstract approach outlined in Remark \ref{th:theory}, this has potential applications to Fukaya categories.
\end{remark}

\section{Additional remarks on sign conventions\label{sec:signs}}
As mentioned briefly in Section \ref{sec:tqft}, we want an $A_\infty$-algebra with vanishing compositions of order $>2$ to be exactly the same as a differential graded algebra (with the standard Koszul signs). The sign conventions therefore follow \cite{gj,k} while differing from \cite{seidel04}. Explicitly, if $\mu^d$ is an $A_\infty$-structure in our sense, translation into an $A_\infty$-structure $\bar{\mu}^d$ in the sense of \cite{seidel04} is achieved by setting
\begin{equation} \label{eq:mu-mu}
\bar\mu^d(a_d,\dots,a_1) = (-1)^{|a_1| + 2|a_2|+\cdots+d|a_d|} \mu^d(a_d,\dots,a_1).
\end{equation}
In Floer theory, signs appear through gluing formulae for determinant lines of (real) elliptic operators. In particular, Koszul signs are natural from that point of view, since the determinant line should be understood as an object whose (formal) parity is given by the index of the elliptic operator. There is an additional source of signs, namely the orientations of various moduli spaces of Riemann surfaces. Our convention for $A_\infty$-structures arises from certain choices of orientations of the Stasheff polyhedra $\scrR^{d+1}$, whereas the convention from \cite{seidel04} is obtained by making an additional ad hoc sign change; one can see this by comparing \eqref{eq:mu-mu} with \cite[Equation (12.24)]{seidel04}.

The same principle applies to other constructions in the paper. For instance, consider the first three terms on the right hand side of \eqref{eq:phi2}, which are all of the same overall form. For the one-dimensional Stasheff moduli space $\scrR^4$ underlying $\mu^3(a_3,a_2,a_1)$, we choose an orientation by varying the marked point corresponding to $a_3$ along the boundary of the disc (with its natural orientation), while keeping the other three marked points fixed. Similarly, for the two-dimensional moduli space $\scrR^{3,1}$ underlying $\phi^2(b,a_2,a_1)$, we vary the interior marked point corresponding to $b$ over the disc (with its natural orientation), while keeping the boundary marked points fixed. As shown in Figure \ref{fig:hex}, the compactification $\bar\scrR^{3,1}$ has three boundary strata which are naturally identified with copies of $\scrR^4$. Figure \ref{fig:halfplane} shows the Riemann surfaces described by points in $\scrR^{3,1}$ close to  each of these boundary strata (for reasons which will become clear presently, it is convenient to draw these surfaces as punctured upper half planes rather than discs; the marked point at ``$+i\infty$'' corresponds to the output of the operation). Deforming the surface towards the boundary stratum corresponds to moving the interior marked point downwards; deforming it parallel to the boundary stratum (in the direction given by the orientation of $\scrR^4$) corresponds to moving the rightmost boundary marked point to the right. One observes that the resulting local coordinates are compatible with the orientation we have chosen for $\scrR^{3,1}$ in the first and third case, and opposite in the second case; this explains the degree-independent part of the signs with which these terms appear in \eqref{eq:phi2}. Figure \ref{fig:halfplane2} indicates the parallel argument for the other three terms on the right hand side of \eqref{eq:phi2}. Formulae involving Dehn twists, such as \eqref{eq:k-kappa}, can be analyzed in the same way (that particular one is actually built from the same underlying moduli space $\scrR^{3,1}$).
\begin{figure}
\begin{centering}
\begin{picture}(0,0)(-50,0)%
\includegraphics{halfplane.pstex}%
\end{picture}%
\setlength{\unitlength}{3158sp}%
\begingroup\makeatletter\ifx\SetFigFont\undefined%
\gdef\SetFigFont#1#2#3#4#5{%
  \reset@font\fontsize{#1}{#2pt}%
  \fontfamily{#3}\fontseries{#4}\fontshape{#5}%
  \selectfont}%
\fi\endgroup%
\begin{picture}(5148,4258)(39,-3718)
\put(-606,139){\makebox(0,0)[lb]{\smash{{\SetFigFont{10}{12}{\rmdefault}{\mddefault}{\updefault}{$\mu^3(a_2,a_1,\phi^0(b))\,:$}%
}}}}
\put(-606,-1539){\makebox(0,0)[lb]{\smash{{\SetFigFont{10}{12}{\rmdefault}{\mddefault}{\updefault}{$\mu^3(a_2,\phi^0(b),a_1)\,:$}%
}}}}
\put(-606,-3299){\makebox(0,0)[lb]{\smash{{\SetFigFont{10}{12}{\rmdefault}{\mddefault}{\updefault}{$\mu^3(\phi^0(b),a_2,a_1)\,:$}%
}}}}
\put(2906,239){\makebox(0,0)[lb]{\smash{{\SetFigFont{10}{12}{\rmdefault}{\mddefault}{\updefault}{\color[rgb]{0,0,0}rescale}%
}}}}
\put(2906,-1411){\makebox(0,0)[lb]{\smash{{\SetFigFont{10}{12}{\rmdefault}{\mddefault}{\updefault}{\color[rgb]{0,0,0}rescale}%
}}}}
\put(4326,-2911){\makebox(0,0)[lb]{\smash{{\SetFigFont{10}{12}{\rmdefault}{\mddefault}{\updefault}{\color[rgb]{0,0,0}towards $\partial\bar{\scrR}^{3,1}$}%
}}}}
\put(4326,-3211){\makebox(0,0)[lb]{\smash{{\SetFigFont{10}{12}{\rmdefault}{\mddefault}{\updefault}{\color[rgb]{0,0,0}parallel to $\partial\bar{\scrR}^{3,1}$}%
}}}}
\end{picture}%
\caption{\label{fig:halfplane}}
\end{centering}
\end{figure}
\begin{figure}
\begin{centering}
\begin{picture}(0,0)%
\includegraphics{halfplane2.pstex}%
\end{picture}%
\setlength{\unitlength}{3158sp}%
\begingroup\makeatletter\ifx\SetFigFont\undefined%
\gdef\SetFigFont#1#2#3#4#5{%
  \reset@font\fontsize{#1}{#2pt}%
  \fontfamily{#3}\fontseries{#4}\fontshape{#5}%
  \selectfont}%
\fi\endgroup%
\begin{picture}(5134,1188)(439,-268)
\put(451,764){\makebox(0,0)[lb]{\smash{{\SetFigFont{10}{12.0}{\rmdefault}{\mddefault}{\updefault}{\color[rgb]{0,0,0}$\phi^1(b,\mu^2(a_2,a_1))\,:$}%
}}}}
\put(2251,764){\makebox(0,0)[lb]{\smash{{\SetFigFont{10}{12.0}{\rmdefault}{\mddefault}{\updefault}{\color[rgb]{0,0,0}$\mu^2(a_2,\phi^1(b,a_1))\,:$}%
}}}}
\put(4051,764){\makebox(0,0)[lb]{\smash{{\SetFigFont{10}{12.0}{\rmdefault}{\mddefault}{\updefault}{\color[rgb]{0,0,0}$\mu^2(\phi^1(b,a_2),a_1))\,:$}%
}}}}
\end{picture}%
\caption{\label{fig:halfplane2}}
\end{centering}
\end{figure}

\vspace{.25 cm}
Department of Mathematics \\
Massachusetts Institute of Technology \\
77 Massachusetts Avenue \\
Cambridge, MA 02139, USA \\

\vspace{.25 cm}
Institute of Mathematics \\
Hebrew University, Givat Ram \\
Jerusalem, 91904, Israel \\

\end{document}